\documentclass[11pt]{article}
\usepackage{amsmath,amsfonts,amssymb,graphicx,amsthm}
\usepackage{epstopdf}
\usepackage{verbatim,bbm,latexsym}
\usepackage{epsfig,color,ifthen}
\usepackage[margin=1in]{geometry}

\newcommand{\bff}[1]{{\mbox{\boldmath$#1$}}}

\newtheorem{remark}{Remark}
\newtheorem{theorem}{Theorem}
\newtheorem{Lem}{Lemma}
\newtheorem{Exa}{Example}
\newtheorem{Pro}{Proposition}
\newtheorem{corollary}{Corollary}
\newtheorem{Def}{Definition}
\def\BEN{\begin{enumerate}}
\def\EEN{\end{enumerate}}   \def\im{\item}
\def\beq{\begin{eqnarray}} \def\eeq{\end{eqnarray}}
\def\le{\left} \def\ri{\right}
\def\te#1{\mathrm{e}^{#1}}
\def\T{\widetilde}  \def\H{\widehat}
\def\bar{\overline}

\def\I{\infty}  \def\a{\alpha}
 \def\b{\beta}

\def\g{\gamma}  \def\d{\delta} \def\de{\delta}
  \def\th{\theta}
\def\e{\epsilon} \def\k{\kappa} \def\l{\lambda}  \def\n{\eta}
\def\Hr{\hat{\rho}}  
  \def\r{\rho}    \def\D{\Delta} \def\K{\mbb K}
\def\t{\tau}   \def\f{\phi}   

\def\rui{\psi}

\def\Rui{\Psi} 
\def\fz{f_C(0)}
\def\Lra{\Longrightarrow}  
\long\def\symbolfootnote[#1]#2{
\begingroup
\def\thefootnote{\fnsymbol{footnote}}\footnote[#1]{#2}
\endgroup}
\def\fn{\symbolfootnote}
\def\fr{\frac}  \def\kil{\mathcal E}
\def\PH{phase-type }  \def\F{\Phi}
\newtheorem{Qu}{Problem}
\def\beQ{\begin{Qu}}
  \def\eeQ{\end{Qu}}  \def\CL{Cram\'er-Lundberg }

\def\P{{\mathbb P}} \def\qq{\quad}
    \def\E{{\mathbb E}} \def\sec{\section}
    \def\ssec{\subsection}\def\la{\label}
    \def\beXa{\begin{Exa}} \def\eeXa{\end{Exa}}
\def\beR{\begin{remark}} \def\eeR{\end{remark}}
\def\beq{\begin{eqnarray}} \def\eeq{\end{eqnarray}}
\def\bea{\begin{eqnarray*}}
\def\eea{\end{eqnarray*}}
\def\le{\left} \def\ri{\right}
\def\te#1{\mathrm{e}^{#1}}
\def\T{\widetilde} \def\Tl{\tilde{\lambda}}  \def\Tq{\T{\q}}
\def\Eq{\Leftrightarrow} \def\beT{\begin{theorem}}
  \def\eeT{\end{theorem}}  \def\beL{\begin{Lem}} \def\eeL{\end{Lem}}
  \def\mT{{\mathcal T}}  
  \def\mJ{{\mathcal J}} \def\und{\underline}
     \def\R{{\mathbb R}}
    \newcommand{\dd}{\mathrm d} \newcommand{\matI}{\mathbb{I}}
    \def\vt{\vartheta} \def\no{\nonumber} \def\ovl{\overline}
    \def\mbb{\mathbb} \def\mbf{\mathbf} 
  \def\q{q} \def\p{r} 
\def\bc{\begin{cases}
  }
\def\ec{\end{cases}}
 
\def\pt{\t_{r}}  \newcommand{\diff}{{\rm d}}
\begin{document}
\title{On the central  management of  risk networks}

\author{Florin Avram\thanks{Laboratoire de Math\'ematiques Appliqu\'ees, Universit\'e de Pau,  France, email: {\tt florin.avram@univ-pau.fr}} \hspace{30pt}Andreea Minca \thanks{Cornell University, School of Operations Research and Information Engineering, Ithaca, NY, 14850, USA, email: {\tt acm299@cornell.edu}}}
\maketitle

\begin{abstract}
This paper provides a methodology for the central management of risk networks, based on the  extensive one-dimensional machinery available  for one-dimensional risk models.
Specifically, we define efficient subsidiaries of a central branch, a concept that takes into account the cumulative dividends generated by the subsidiaries, as well as the cost of bailing them out by the central branch.
In the case of deterministic central branches with a single subsidiary, we find closed form solutions for the value of the subsidiary to the central branch.  Moreover, we extend these results in the case of hierarchical networks.
In the case of non-deterministic central branches with one subsidiary, we compute  approximate value functions
 by applying  rational approximations,  and by using the recently developed  matrix scale  methodology.

{\bf AMS 2000 Math. Subject Classification:} Primary 60G51, Secondary 60K30, 60J75
\end{abstract}
\tableofcontents

\section{Introduction}

{\bf Motivation}.
In recent years considerable effort has been devoted in the finance and insurance literature to model default (or ruin) in the context of a system of interacting firms.
The insurance giant AIG was almost led to  bankruptcy after its subsidiary AIG Financial Products faced contingent claims on its derivatives positions and forced a massive government bailout in 2008. AIG Financials Products, liquidated in 2011, will probably remain as a textbook example of a nonviable subsidiary whose losses were  disproportionate with respect to its prior profitability.

Similar  to the problem of a central branch in need to bail out its subsidiaries is the problem of central clearing counterparty, mandated by the Dodd-Frank regulations in response to the recent crisis. Just as central branch can be ruined by one subsidiary, a central clearinghouse can  default if its capital is insufficient to cover the default of one (or several) members. This in itself could be a systemic event, so it is a current concern if central counterparties are effectively managed.

Another example which falls under the umbrella of a central branch is a reinsurance company which could face potentially large claims, depending on the primary insurers' deficit.

In this paper we consider a general model of central branch risk networks.
A key point is  defining a notion of efficiency of a subsidiary. Whereas intuitively it is clear that not all subsidiaries represent viable businesses from the perspective of the central branch, and some should be closed, it is far from clear which criteria should be utilized to achieve that. We will illustrate this by  comparing several such criteria.

Let us turn now to the mathematical framework of our  study.

{\bf   Multi-dimensional risk networks } (\textbf{MRN})  are defined by:
 \begin{eqnarray*} &&\bff X(t)=\bff u + \bff
c t- \bff S(t)=(X_i(t),\ t\geq 0, i \in {\mathcal{I}}), \end{eqnarray*}
where ${\mathcal{I}} = [1, ..., I]$ is a finite set,
  the vector $\bff u$ represents the
capital of the \textbf{MRN} \,at time $0$,   the vector $\bff c$
represents a constant cash inflow rate,
and  $\bff S(t)$ is a 
process representing  cash outflows at time $t$, which may include both  Levy and  Sparre-Andersen renewal components.

If no boundary condition is specified, we will call this a  {\bf free spectrally negative  \textbf{MRN}}.

\begin{remark} The minus sign comes from the  one-dimensional case  most  studied historically,   the spectrally negative Cram\'er-Lundberg process, but the case when $X(t)$ is spectrally positive is also interesting. The case of spectrally two-sided $X(t)$ is of course  interesting, but harder. \end{remark}

\begin{remark} The simplest case  is    when  the component are i.i.d. compound renewal Sparre-Andersen processes, generated by i.i.d. pairs
 of inter-arrival times and claims $(A^{(i)}_j, C^{(i)}_j), j=1,2, ...$
\begin{eqnarray*}\label{SA}
X_i(t)= u_i + { c_i} \; t -S_i(t), \;
S_i(t)=\sum_{j=1}^{N_i(t)} C^{(i)}_j,\; \; i =1, ....,I, \end{eqnarray*} where $N_i(t)=\max \{k: T^{(i)}_k:=\sum_{j=1}^k A^{(i)}_j \leq t\}$ are     renewal counting  processes  associated to the  independent inter-arrivals (see e.g. \cite{AA}), with intensity $\lambda_i=E (A^{(i)}_1)^{-1}, $ and the claims sizes $C^{(i)}_j, j \geq 1$ are
nonnegative i.i.d. random variables with arbitrary marginal
distribution functions denoted by $F_{i}(x), i \in {\mathcal{I}}$, with finite expectation, denoted $m_i$. However, the restriction to jump processes is not essential.

\end{remark}
\begin{remark}  After reaching special subsets exterior to the
state space, several  continuation/regulation mechanisms are possible, like absorption, reflection, or jumping to the interior. These correspond to various
possible interactions between the components at times of distress. \end{remark}

\begin{Exa}{\bf A toy  example  with one absorbing boundary and several reflecting boundaries}. Consider   a central branch which   must simultaneously manage several
 subsidiaries.

 The central branch will  keep the
 subsidiaries  solvent by bail-outs until the moment of its bankruptcy, or according to a pre-specified rule, e.g. until the subsidiary is deemed non-profitable.
The subsidiaries will pay dividends to the central branch. Finally, the expected present (net) value  to the central branch of a subsidiary consists
 in the difference between the expected discounted payments and  bail-out amounts.

\end{Exa}

 This example
  suggests the following  model:
 \begin{Def}{\bf A
central  branch (CB) network \label{ex:cb}} is formed from:
\BEN \im  A unit,
called  central branch, with reserves denoted by
 $X_0(t)$, whose ruin time
$$
\tau_0 = \inf\{t\ge 0: X_0(t) < 0\}
$$
causes the ruin of the whole network.

\im  Several
subsidiaries  $X_i(t), i=1,\ldots,I$ 
that must 
be kept nonnegative or above certain prescribed levels, by  transfers from the CB.%
\EEN
\end{Def}
For this network, the boundaries $u_i=0, i=1,...,I$ are reflecting  and $u_0=0$ is absorbing.

\begin{remark} There are many applications of the central branch concept:  a government/central bank, a reinsurance company, an insurance group, a central clearinghouse, etc.

Another  interesting application is that
 of a coalition or default fund  created by    several institutions,   for bailing them out  when they go bankrupt.
 Interesting issues here are determining fair conditions for merging into (profit participation schemes) and  splitting out of  the coalition.

One of the most important activities of central branches  is acquiring and closing subsidiaries;
A critical building block in understanding the management of the CB risk network is modeling one subsidiary managed by the CB. The benefit is represented by the dividends received from the subsidiary, while the costs are the cumulated bail-out costs, the number of bailouts, up to a certain time set by the CB and referred to as ``patience".
Establishing decoupled rules for the management of subsidiaries is the key step towards solving  problems pertaining to the whole risk network, such as  maximizing  the total cumulated net income from all the the subsidiaries, until the eventual bankruptcy of the CB.

\end{remark}

{\bf Notation}. Denote the ruin times and ruin probabilities (finite time and eventual)
of the components when isolated from the network by\begin{eqnarray*}&&\t_i(u_i) = \inf\{t\ge0: X_i(t) < 0\}, i=1,2, ...\\&&\Rui_i(t,  u_i) = P(\tau_i(u_i) < t), \; \; \;
\Rui_i(u_i):= P(\t_i(u_i) < \I).  \label{onedim}\end{eqnarray*}
The ruin probability of the CB and  its Laplace transform will be  denoted respectively  by
 \beq &&\Rui(t, \bff u)=\Rui(t, \bff u, \bff c)= P_{\bff u}[\t_0 < t], \;
\H{\rui}_\q(  \bff u):=\H{\rui}_\q(  \bff u, \bff c) =E_{\bff u}\left[ e^{-\q \t_0}\right]. \label{def_ruinLT} \nonumber \eeq

  {\bf One-dimensional first passage problems} have been very extensively studied; typically, Laplace transforms are available, especially when either $A_i$ or $C_i$ have a matrix exponential distribution, and  explicit inversion of the Laplace transforms is  also possible sometimes, especially when at least one of $A_i$ or $C_i$ have an exponential distribution.

 {\bf The eight pillars of one-dimensional first passage problems for spectrally negative (and positive) processes}.  In the last decade it was realized that the solutions of a large gamut   of  problems  (dividends, drawdowns,  barrier options, exotic options, etc...)  may be expressed in terms of the solutions of eight basic problems \eqref{onesided}, \eqref{s}, \eqref{excscale}, \eqref{refbailouts}, \eqref{div}, \eqref{sevruin}, \eqref{regsevruin}, \eqref{resolvent_density}, which at their turn  may be   ergonomically  expressed in terms
 of a couple of basic solutions called  {\bf scale functions}\fn[4]{The eight pillars are respectively  the solutions of the one, two-sided and drawdown constrained smooth first passage problems, the formulas  of reflected bailouts, dividends, severity of ruin and regulated severity of ruin, and the resolvent density.}.

  These eight  pillars expand   the content of the ``basic dictionary entry" for a specific Markovian stochastic process furnished by its  generator, or by its equivalents: the symbol and the first scale function, and act as a sort of a instruction kit for  solving a variety of  first passage problems.

  The precursor of this  idea may be found  in \cite{AKP}, where the "second scale function" was introduced.  Its full   development may be found  in
\cite{pistorius2005potential}, where  a list of identities expressible in terms of the two scale functions was offered --see also \cite{KKR}. Remarkable  extensions which involve ``relaxed" first passage times   were provided in \cite{AIpower,AIZ}.   It seemed  probable from the start that these identities could be applied to  Markov processes  with state dependent generators  as well, but the difficulty to compute the scale functions stopped further research, with  the notable exception of \cite{albrecher2010algebraic}.

The idea of constructing solutions out of basic scale functions came finally to full bloom while being  extended  to spectrally negative Markov additive processes (which include L\'evy processes in a modulated Markovian  environment, and  the case of Sparre-Andersen processes with phase-type arrivals) in \cite{KP,Iva,IP},
  reviewed in the appendix. This extension requires {\bf matrix scale functions}, and the  end result is a Mathematica {\bf matrix scale package}
  offered by \cite{Ivapackage}.

   {\bf  From one to several dimensions}.  Multi-dimensional first passage problems are   considerably  harder than   { one dimensional ones}, and
 one cannot  expect  general formulas\fn[4]
{One  exception is a
  Pollaczek-Khinchine type formula for the transform of ruin probabilities  $\Rui(\bff u)$
of spectrally negative networks  provided in the foundational paper
\cite{CYZ}. However, this formula involves several  unknown functions (the Laplace transforms over each boundary  facet of the state space), and  it isn't at all  obvious how to  exploit this
  formula numerically.}. In this paper we provide   approximations in the multidimensional case by reduction to  {\bf   one-dimensional results}.

  {\bf Efficient subsidiaries}.
A crucial issue for a coalition is how to accept {\it efficient members} and eventually reject them if they are not  efficient.
For that it is natural to evaluate each member separately, by classic one dimensional risk measures like ruin probabilities, or  the value of  future dividend payments made to the coalition.
However, the choice of an economic principle for evaluating efficiency is not at all obvious.

One often used  approach is to optimize  discounted dividends  of rate
 $d =c \g, \; \g \in (0,1] $ taken above a constant threshold $b$ ($\g$ represents the proportion of income taken above the threshold) -- see for example \cite{albrecher2007exact,avanzi2009strategies}. The process obtained by subtracting the dividends is called reflected when $\g=1$,   and refracted when $\g \in (0,1) $.

  {\bf Judging efficiency as readiness to pay dividends}.  Another
   method to judge efficiency  was made in the conference paper \cite{AM15},  who proposes    to
 define efficiency   as {\bf local optimality of $b=0$} over some interval  $[0,\e), \e>0.$

 The motivation is that such subsidiaries are functional from the start and can contribute  cash-flows to the central branch without having to wait first until its reserves build out; effectiveness is thus translated in this paper as {\bf readiness}.  A second motivation is that inefficient subsidiaries may be turned into efficient by setting up rules  to monitor their time passed in
 in "orange zones", and eventually close them when necessary. This may be achieved by ``killing" them with a rate $\th_i$ in the orange zone, where $\th_i$ is chosen to render  the barrier $b=0$
 locally  optimal.  Furthermore, one may use "two step" killing rates, and choose the killing rate as $\I$ in "red zones", as suggested in \cite{czarna2015note}. We show in a variety of example that this procedure produces reasonable results.

 Summing up, we will propose that a subsidiary
 is: \BEN   \im  {\bf Non-efficient} and rejected immediately if its loading factor $\fr{c }{\l E[C_1]}-1:=\r^{-1} -1$ is not nonnegative, since this  implies an infinite number of bail-outs.
 \im {\bf Totally efficient} and accepted for ever in the coalition iff \beq \la{e:efd} \r=\l E[C_1]/c <1,  \text{  and } k \leq f(\q) \eeq
 where $\q$ is the discount rate, and $f(\q) $ is an increasing function of $\q$, with $f(\q) <0,$ obtained as optimality of $b=0$ for some specific dividends distribution scheme, and $k \geq 1$ captures the cost associated with capital infusions towards a subsidiary.
 \im {\bf Partially efficient} if the  loading condition $\r <1$ is  satisfied, and $k > f(\q) .$
  These subsidiaries will also be accepted, but only with an {\it impatience rate} $\th,$ resulting in killing the subsidiary after its time or lowest value  in an orange zone exceeds an  exponential r.v.  of rate
 $\th$.  The impatience  rate $\th$  is  chosen so that
 $$ \T f(\d,\th)=k,$$ where $ \T f(\d,\th)$ is computed from the optimality of $b=0$ for the {\it impatience modified} value function.
 This is illustrated in example \ref{e:th} below, where $ \T f(\d,\th) =f(\q + \th)$.
\EEN

  {\bf
  The classic De Finetti objective}  maximizes expected discounted dividends until the ruin time.   For judging efficiency, it is natural to take also into account  a final bail-out, resulting in the optimization objective:
   \beq\label{e:robj} V^{(w)}(x)=\sup_{\pi} E_x \le[ \int_0^\t e^{- \q t} d D^{\pi}(t) +  e^{- \q \t} w(U(\t))\ri], \eeq
   where $w(u)$ is a so called Gerber-Shiu penalty function.

   The optimal dividend distribution  is  of {\it multi-barrier} type \cite{gerber1972games}, and the  end result may  be expressed in  terms of  scale functions \cite{APP,APP15}.  Further conditions are necessary to ensure that  {\bf single constant barrier}  strategies suffice \cite{AM,APP,loeffen2008optimality,LR}.

   {\bf Judging efficiency of subsidiaries  over an infinite horizon by optimizing bail-outs and dividends}. Over  an infinite horizon, the subsidiary will possibly need to be bailed out  a number  $N_B \geq 0$ of times.

 In the case of linear transaction costs $k u -K$, the optimization objective   (of particular interest in a bail-out setting) becomes the expectation over an infinite horizon of a linear combination of  discounted dividends $D(t)$, cumulative bailouts $Z(t)$,  and number of bailouts $N_B^{\pi}(t)$  up to time $t$:
  \beq\label{e:Sobj} V^{(k)}(x)=\sup_{\pi} E_x \le[ \int_0^\infty e^{- \q t} d D^{\pi}(t) -k \int_0^\infty e^{- \q t} d Z^{\pi}(t) - K \int_0^\infty e^{- \q t} d N_B^{\pi}(t) \ri].  \eeq

  Since in a diffusion setting this objective has first been considered by Shreve, Lehoczky,  and Gaver (SLG) \cite{shreve1984optimal} -- see also Lokka and
Zervos \cite{lokka2008optimal} -- we will call it the SLG objective.

 For spectrally negative Levy processes, the optimal dividend distribution for  the SLG objective is always of {\bf constant barrier} type, and the  end result may  be expressed in  terms of  scale functions \cite{APP}.

 {\bf Bail-out intervention  times}.  For  bail-out  times, one may consider the classic ruin time $\t=\t_0^-$, and  also  several interesting alternatives generalizing it:
 \BEN

 \im One may replace $\t$ in \eqref{e:robj}  by a  ruin time $\t_\p$  observed with Poissonian frequency $\p$-- see for example \cite{AIZ,albrecher2015strikingly}, which is equivalent to a  Parisian ruin time with exponential grace period -- see for example \cite{LRZ-0}.

   \im  One may further replace $\t$ in \eqref{e:robj}  by a
   {\it bankruptcy time} $\t_{\p,a}$ involving Poissonian-Parisian   grace period below $0$ and also absolute killing at $a<0$ -- see for example \cite[Sec. 4]{Ren}.

\EEN

\noindent {\bf  Some definitions  of efficiency based on prior literature}.

\begin{Exa}{\bf A simple, but unsatisfactory  definition of  efficiency}. Instead of our notion of {\it readiness}, consider defining efficiency as nonnegativity of the SLG objective for the $0$ barrier, when dividends are $D(t)= d \; t$. By \cite[Thm. 1, (4.4)]{APP}
\beq V^{(k)}(0)=\fr c {q}- k \fr{\l m_1} {q}  \geq 0 \eeq
yields $k \leq  \fr c {\l m_1} = \r^{-1}. $
This generalizes easily for the $0$ threshold,
yielding
$$k \leq  \fr d {\l m_1} = \r^{-1} \g.$$
\end{Exa}

\begin{Exa}\label{e:th} {\bf SLG readiness}. Consider now  efficiency defined as the optimality of $b=0$
for the SLG objective  with {\bf constant reflecting barrier}.
 This  problem is  fully analyzed in \cite[ Thm. 3]{APP}, and in particular
 \cite[Lem. 2]{APP}  shows that  the optimal SLG constant barrier  is  $b^*=0$
iff
$  k \leq 1+ \fr {q}{\l}.$

By this criterion,
a subsidiary $i$ with \beq\label{e:Scb} k_i \leq 1+ \fr {q}{\l_i}\eeq will be deemed totally efficient and accepted for ever in the coalition\fn[4]{Otherwise \cite[(5.6)]{APP} \begin{equation}
b^* = \inf\{a>0:[k Z^{(q)}(a) - 1]W^{(q)\prime}(a) - k  \q W^{(q)}(a)^2\leq 0\}  >0\label{dstar}
\end{equation}
 (i.e. the subsidiary must be forgiven from paying dividends in between $(0,b^*)$. See also \cite{kulenko2008optimal,eisenberg2011minimising} for related results.}. Subsidiaries with loading condition $\r_i <1$ and
 $$ k_i > 1+\fr {\q}{\l_i}$$ will be deemed {\bf partially efficient}  and
  accepted
 only for a  random time with law $\kil(\th_i)$,
where  \beq\label{e:bal} \th_i +  {\q}=\l_i ({k_i - 1})\eeq
(rendering thus $b_i^*=0$ optimal with respect to the
 total discount rate $\d_i = \q + \th_i$).

Unfortunately, the criterion \eqref{e:Scb} does not take into account the claim size law.

\end{Exa}

To avoid these shortcomings, we introduce a new  efficiency concept in Section \ref{s:eDeF}.

{\bf Contents and contributions}.
 Our first contribution is to
provide in Section \ref{s:eff} a definition  of  efficiency which is acceptable for a wide variety of subsidiaries satisfying the conditions of Lemma 1,  which include exponential claims.

The central branch model  is described in more detail in Section \ref{s:CB}. Our second contribution, Theorem  \ref{easy} in Section \ref{s:red}, applies to   a {\bf purely  deterministic  CB} $\T X_0(t)=u_0 + c_0 t$ with  one subsidiary. In this case,  the computation of the finite time ruin probabilities and other performance measures (including total subsidiary dividends until ruin)  reduces   to the corresponding computation for a subsidiary with modified initial capital $u_0+ u_1/k_1$ and initial income rate  $c_0+ c_1/k_1$\fn[4]{For other potentially  useful explicit computations  see     \cite{APP08,badescu2011two}.}.  More precisely,
\beq\label{e:Vred} &&\Rui(t,\bff u, \bff c) = \Rui_1(t,u_0/k +
u_1,c_0/k + c_1), \nonumber \\&& V(\bff u, \bff c) = V_1(u_0/k +
u_1,c_0/k + c_1),\eeq
{\bf without any distributional assumptions}!  For example, the subsidiaries may be dual risk processes, or spectrally two-sided Levy processes,...

The   proof of this result, via a   pathwise argument, yields also      an  upper bound when $I>1$, and
an extension to hierarchical networks is given in  Corollary \ref{c:easy}. The notion of hierarchical networks has been used ever since Gerber (\cite{gerber1984}) to model chains of reinsurers and is the most important application of the deterministic CB case.

 Our third contribution, in Section \ref{s:rho}, deals with {\bf non-deterministic  CB}'s with one subsidiary, which do not admit an exact solution due to their complex dependent Sparre-Andersen structure.
We isolate  the CB and one subsidiary, and allow for a non-deterministic component of the CB process to represent  other claims of the CB, for example aggregate  net flows of the other subsidiaries or other liquidity needs of the CB.
 We propose an approximation approach
 based on the classic idea of rational approximations to replace the Sparre-Andersen structure by  a  Markov  modulated  L\'evy structure,  and by using subsequently the matrix scale  methodology.

 We propose  to obtain  {\bf SNMAP} (spectrally negative Markov additive process)  approximations
    for  non-MAP   central branches, by using {\bf bivariate phase-type}  approximations  for the {\bf joint law of the downward ladder time and height}.
The advantage of this approach  is that  once a SNMAP approximation
    is obtained, many similar problems may be solved just by applying the scale
    matrix methodology developed by \cite{KP,Iva,IP}, and using the SNMAP Mathematica  package of J. Ivanovs\cite{Ivapackage}.  Different problems
    are thus solved simultaneously!

  A numeric illustration is performed  in section \ref{s:PDV}, where we consider the problem of choosing a barrier $B$ maximizing CB dividends until ruin, in the  case of one subsidiary with exponential claims (and without dividends).

    In this case,
 {\bf univariate phase-type}  approximations   of the downward ladder density, obtained via a continued fraction expansion 
 --see section \ref{s:rho}\fn[4]{Note that even though here the  CB claim arrival times
 have  an explicit Bessel-type law, 
  our approach replaces this by a phase-type approximation.}, provide a  SNMAP  approximation of the CB process.
Subsequently,  using the SNMAP package of Ivanovs provides the optimal barrier.

  Our methodology based on the three ideas proposed above lays the grounds for an  approximate optimization of risk networks. The first step consists of  setting  impatience parameters $\th_i$ (but not the  initial allocation capital infusion from the CB $v_i$ and dividend barrier $b_i$) for each subsidiary viewed in isolation from the network.
The second step will consist in setting  dividend barrier  parameters $b_i$ and  optimal allocation parameters $v_i$, using a  {\it decoupled} objective of the form
 $V(\bff v)=\sum_{i} V_i(u_i + v_i)$, representing the sum of  benefits to the network from all components, under the constraint $\sum_i v_i =u_0,$ keeping the $\th_i$, and using  the barriers $b_i^*=0$ from the first step as initial  values for an iterative procedure.    When transaction costs are present,  they may   also be incorporated via  the reduction result \eqref{e:Vred}.
A  harder problem would be  to take into account   the possible  bankruptcy of the central branch.

  Other alternatives of  efficiency  based on threshold strategies,
linear reflecting barriers, Poisson evaluation periods, and killing   based on total current bail-out
will be  investigated in future work. To prepare this, a review of the SNMAP matrix scale  approach is provided in the last three sections.

\section{Preliminaries on first passage}\label{sec:MAPexit}

Consider a spectrally-negative MAP $(X(t),J(t))$ and assume that none of the underlying L\'evy processes $X^i_{t}$ is a.s. non-increasing (see Section \ref{s:MAP} for a review of these processes).
Define the first passage times by
\begin{align*}
 &\tau_b^+=\inf\{t\geq 0: X(t)>b\},&\tau_b^-=\inf\{t\geq 0: X(t)<b\},
\end{align*}
 with $\inf\emptyset=+\infty$. We will sometimes write $\t$ for the "ruin time" $\t_0^-$.

 Let $\t_{a,b} = \t_{a,b}= \t_{a}^- \wedge \t_b^+$
denote the ``two-sided'' exit time from an interval $[a,b]$. 
\bigskip

If the process is only observed   at the arrival times $\mT_\p=\{T_i,i=1,2,...,$  of an independent Poisson process of rate $\p$,   the analog concepts are the stopping times
\begin{eqnarray*}T_b^+= \inf\{T_i:\; X(T_i)  > b\}, \quad T_{b}^- = \inf\{t > 0:\; X(T_i) < {b}\}
 \end{eqnarray*}
  see for example \cite{AIZ,albrecher2015strikingly}.

 A related concept is the {\bf Parisian ruin time}
 $ \pt(a)$ defined as the first  time when the most recent excursion  below $a$
 has exceeded an exponential rv $\kil_\p$ of rate $\p$ (in other words, the process is killed below $a$ at rate $\p$) -- see for example \cite{LRZ-0}. When $a=0$, the notation $\pt$ will be used.

{\bf
  The classic De Finetti objective}  maximizing expected discounted dividends until the ruin time may be taken as a basic principle for judging efficiency. In this context, it is natural to take also into account  a final bail-out, resulting in the optimization objective:
   \beq\label{e:robj} V^{(w)}(x)=\sup_{\pi} E_x \le[ \int_0^\t e^{- \q t} d D^{\pi}(t) +  e^{- \q \t} w(U(\t))\ri], \eeq
   where $w(u)$ is a so called Gerber-Shiu penalty function.

   The optimal dividend distribution  is  of {\it multi-barrier} type \cite{gerber1972games}, and the  end result may  be expressed in  terms of  scale functions \cite{APP,APP15}.  Further conditions are necessary to ensure that  {\bf single constant barrier}  strategies suffice \cite{AM,APP,loeffen2008optimality,LR}.

    Over  an infinite horizon, the subsidiary will possibly need to be bailed out  a number  $N_B \geq 0$ of times.

 In the case of linear transaction costs $k u -K$, the optimization objective   (of particular interest in a bail-out setting) becomes the expectation over an infinite horizon of a linear combination of  discounted dividends $D(t)$, cumulative bailouts $Z(t)$,  and number of bailouts $N_B^{\pi}(t)$  up to time $t$:
  \beq\label{e:Sobj} V^{(k)}(x)=\sup_{\pi} E_x \le[ \int_0^\infty e^{- \q t} d D^{\pi}(t) -k \int_0^\infty e^{- \q t} d Z^{\pi}(t) - K \int_0^\infty e^{- \q t} d N_B^{\pi}(t) \ri].  \eeq

  Since in a diffusion setting this objective has first been considered by Shreve, Lehoczky,  and Gaver (SLG) \cite{shreve1984optimal} -- see also Lokka and
Zervos \cite{lokka2008optimal} -- we will call it the SLG objective.

 For spectrally negative Levy processes, the optimal dividend distribution for  the SLG objective is always of {\bf constant barrier} type, and the  end result may  be expressed in  terms of  scale functions \cite{APP}.

 {\bf Bail-out intervention  times}.  For  bail-out  times, one may consider the classic ruin time $\t=\t_0^-$, and  also  its "soft" Poissonian-Parisian alternatives generalizations.

\sec{Efficient subsidiaries \label{s:eff}}

As argued in the introduction, we will consider that a company is:
\BEN
\im {\bf Totally efficient} and accepted for ever in the coalition if  \beq \la{e:efd} \r=\l E[C_1]/c <1,  \text{  and } k \leq f(\q) \eeq
 where $\q$ is the discount rate, and $f(\q) $ is an increasing function of $\q$, with $f(\q) <0,$ obtained as optimality of $b=0$ for some specific dividends distribution scheme, and $k \geq 1$ captures the cost associated with capital infusions towards a subsidiary.
 \im {\bf Partially efficient} if the  loading condition $\r <1$ is  satisfied, and $k > f(\q) .$
  These subsidiaries will also be accepted, but only with an {\it impatience rate} $\p,$ resulting in killing the subsidiary after its time or lowest value  in an orange zone exceeds an  exponential random variable  of rate
 $\p$.  The impatience  rate $\p$  is  chosen so that
 $$ \T f(\q,\p)=k,$$ where $ \T f(\q,\p)$ is computed from the optimality of $b=0$ for the {\it impatience modified} value function.
In example \ref{e:th} we have $f(\q) =  1+ \fr {\q}{\l_i}$ (see \eqref{e:Scb}) and we set $ \T f(\q,\p) =f(\q + \p)$.
Motivated by the fact that the criterion \eqref{e:Scb} does not take into account the claim size law, we introduce in
Section \ref{s:eDeF}  an efficiency concept which does not have this shortcoming.
\EEN

We list now  several  alternatives of  dividend payment  strategies of a subsidiary towards a   coalition, which could turn out  to provide  appropriate definitions of readiness  in our context, and which we plan to investigate in the future. 

\BEN \im   {\bf Linear reflecting barrier}
   strategies  consist in specifying a function $b(t)=b + (c-d) t$, and taking as dividends all surpluses above it.

\im {\bf Two-step premium/refraction/threshold dividend} strategies
 postulate that only a proportion $0< \g= \fr d c \leq 1 $ of the premium income may be paid as dividends above the {\bf threshold} $b$ (the name barrier  is replaced  by threshold, since the process goes on evolving above $b$).
When $\g=1,$ the  two-step premium policy reduces to a reflecting  barrier policy.

Threshold  strategies are  motivated by the fact that the  optimal dividend distribution under the De Finetti objective with an extra constraint $$d D(t) \leq d \; dt,$$
where $0< d < c$, is known to be consist in modifying the premium from$c$ to $d$ above a constant barrier $b$ -- see for example \cite{albrecher2009optimality,azcue2012optimal}. Another motivation is provided in \cite{shen2013alternative}.

  For \CL processes, ruin is sure or not depending on whether \cite{boxma2011threshold}
\beq 1-\g < \r:= \fr{\l E [C_1]}{c} \eeq  or not,  where $\l,c$ denote the rate of arrivals and premium, respectively.

For threshold/refraction policies, some   formulas  expressed in terms of  scale functions are provided in \cite{kyprianou2010refracted,kyprianou2013gerber}.

\im A third possibility is that the subsidiary pays "tax" -- a proportion $\g$ of the premium income, whenever the subsidiary is at a running maximum--  see for example  \cite{albrecher2008levy,AACI,AIpower}.
(when $\g=1,$ this reduces again to paying dividends above  a constant reflecting barrier).

\EEN

We examine now one possible efficiency criterion.

 \ssec{Efficiency based on  De Finetti constant dividend barrier and linear penalty \label{s:eDeF}}

As noted above, a critical step in defining efficiency is the choice of the value function. Our choice is based on the case of
De Finetti constant  barrier and linear penalty (with $k$ proportional cost and $K$ the fixed cost), for which the value function is known explicitly, see \cite[(13)]{avram2009optimal}, \cite{LR,APP15}.

 Putting $\bar F(x)=\int_0^x F(u)  du, F(u)=k\le(Z(u)- \k'(0) W(u)\ri)-K \q W(u)$\fn[4]{$F(u)=-K Z^{(0)}(u)+ k Z^{(1)}(u)$, where the coefficients of $-K,k$ in $F(x)$ are found by differentiating  the second scale function \cite{APP15} $ Z(x,\th)= e^{\th x} \le( 1- \k(\th) \int_0^x e^{- \th y} W(y) dy  \ri)$, $0$ and $1$ times respectively, with respect to $\th$, and then differentiating  once more with respect to $x$.},    they find that   the value function
may be written as
 $$ V(x) =\bar F (x) + W(x) G(b), \; G(b)=\fr{1- F(b)}{W'(b)}=\fr{1- k\le(Z(b)- \k'(0) W(b)\ri)+
K \q W(b)}{W'(b)}.$$

 The optimality condition is obtained by differentiating
 the "barrier influence function"  $G(b)$  \cite[Sec. 13.1]{APP15}.
 Putting $ \Tq:=\fr {q}{c}$, $\D(b)=c \le( W^{'2}(b)-W(b) W''(b)\ri) $, we find
 \beq && H(b):=G'(b) {(W'(b))^2}=
 \le(k \le(\k'(0) W'(b)- \q W(b))\ri)+
K \q W'(b)\ri) W'(b) \no
\\&&-W''(b) \le(1+ k (\k'(0) W(b)-1 - \q \bar W(b)) + K \q W(b) \ri)\no \\&&=\le(K \Tq +k (1-\r)\ri) c \le(W^{'2}(b)-W(b) W''(b)\ri) +
k  \q \le(W''(b)\bar W(b)- W'(b)W(b)\ri)   +(k-1) W''(b)
\no
\\&&=k  \le((1-\r) \D(b)+\q \le(W''(b)\bar W(b)- W'(b)W(b)\ri) + W''(b)\ri) +  K \Tq \D(b)   - W''(b). \end{eqnarray}

Now to compute  $ W''(0)$, let $K(s)= \fr 1 c \fr{\T \l s \H{\bar{F}}(s)}{1-\T \l  \H{\bar{F}}(s) -\d/s}, \Tl:=\fr {\l} c$
the expression used in \cite[pg. 33]{kuznetsov2013theory} to find $W'(0)=\lim_{s \to \I} K(s)=\fr{\Tl + \Tq}{c}$.  Following the same approach, we find
\beq  \label{e:secder}  W''(0)=\lim_{s \to \I} s(K(s) -W'(0))= \fr 1 c \le( (\Tl + \Tq)^2 - \Tl \fz\ri).\eeq

Using furthermore $\bar W(0)=0, W(0)=1/c, c \D(0)= \Tl \fz$,  we find:

\bea&&  c G'(0) {(W'(0))^2}=k  \le((1-\r) \Tl \fz -\q  W'(0) + c W''(0)\ri) +  K \Tq \Tl \fz   - c W''(0)\\&&=k  \le(\Tl (\Tl + \Tq) -\r \Tl \fz \ri) +  K \Tq \Tl \fz   - \le((\Tl + \Tq)^2  - \Tl \fz\ri)
\eea

The optimality condition is:
 \beq \label{e:ef} && G'(0) \leq 0 \Eq
 k \leq
 f(\Tq,K):=\Tl^{-1} \fr{(\Tl + \Tq)^2 - \Tl \fz-K \Tq  \Tl \fz}{ \Tq + \Tl (1- m_1  \fz)} \label{e:ef2}
  \eeq

\begin{remark} Note that if $K=0, \fz=0,$ this reduces to SLG readiness. \end{remark}
 \beXa As a check, with exponential claims, the  scale function is:
$$
W^{(q)}(x) = A_+ \te{\zeta^+(q)x} - A_- \te{\zeta^-(q)x},
$$
where $A_\pm = c^{-1} \frac{\mu +
\zeta^\pm(q)}{\zeta^+(q)-\zeta^-(q)},$ and $\zeta^+(q)=\F(q),$
$\zeta^-(q)$ are the largest and smallest roots of the polynomial
$c^{-1}(\psi(s)-q)(s + \mu)=  s^2 - s(\Tl+ \Tq -\mu) - \Tq \mu$:
$$\zeta^\pm(q)=\frac{{\Tq} + \Tl -\mu  \pm
\sqrt{({\Tq} + \Tl-\mu )^2 + 4  {\Tq}\mu}}{2}.$$

Thus, \begin{eqnarray*}&&W''(0)=A_+ (\zeta^+(q))^2 - A_- (\zeta^-(q))^2=\Tq \mu (A_+ - A_-)+  (\Tl+ \Tq -\mu) \le(A_+ \zeta^+(q) - A_-  \zeta^-(q)\ri) \end{eqnarray*}and
\begin{eqnarray*}&&c W''(0)=\Tq \mu +  (\Tl+ \Tq -\mu)(\Tl+ \Tq ) =  (\Tl+ \Tq )^2-\Tl \mu,  \end{eqnarray*}confirming \eqref{e:secder}.

  Finally, with $K=0$ and exponential  claims of rate $\mu,$ we recover  \cite[Lem 13.2 b]{APP15}:
 \beq  k \leq \fr{(\Tl + \Tq)^2 - \Tl \mu}
{\Tl \Tq}=\fr{(\l + \q) ^2 - \l c \mu}{\l \q}  \label{e:exef}\eeq
\eeXa

\beL  With  general  claims admitting a first moment such that $\r \leq m_1 \fz \leq 1$, the function $f(\Tq,K), \Tq >0$ defined in \eqref{e:ef} is increasing  for any $0 \leq K \leq \fr{\r^{-1} -1}{1- m_1  \fz } + \fr 1 {m_1\fz}$. \eeL

{\bf Proof:}  The function defined in the right hand side of \eqref{e:ef} is eventually increasing when $\q \to \I$. Its
value in $0$ is proportional to $\fr{ \r- m_1 \fz}{1-m_1 \fz}$ and hence negative under the hypothesis.   The derivative  with respect to $\Tq$ of
of \eqref{e:ef2} is proportional to
\begin{eqnarray*}(z+1-\fz  m_1)^2- \fz^2  m_1^2 +\fz m_1 \r^{-1}
   +K \fz
   \left(\fz
   m_1- 1 \right), z=\Tq/\Tl.\end{eqnarray*}

   Finally, assuming  $\r \leq 1$, we find three different  conditions ensuring that  the derivative is positive for any $z \geq 0$:
   \bea \bc  K/m_1 \leq \fr{\r^{-1} -1}{1- m_1  \fz } + \fr 1 {m_1\fz}  &\fz  m_1  < 1\\ K \geq 0  &\fz  m_1 = 1\\  K/ m_1 \geq \fr{ (\fz m_1 - \r^{-1})_+}{m_1 \fz  - 1}
   &\fz  m_1 \geq  1\ec \eea

Moreover, in the case $\fz  m_1 \geq  1$ if we have $K/ m_1 < \fr{ (\fz m_1 - \r^{-1})_+}{m_1 \fz  - 1}$ then there are $r_1(K) \leq r_2(K)$ such that the function
$\Tq \to f(\Tq,K)$ defined in \eqref{e:ef} is increasing on $[0, r_1]$, decreasing in $[r_1.r_2]$ and increasing to $\infty$ on $[r_2, \infty)$. Because the image of $\Tq \to f(\Tq,K)$ is $[0, \infty)$ we can always find an appropriate killing rate, e.g. the smallest solution to $k = f(\Tq,K)$ (if $k \in Im_f([r_1, r_2])$, there are three such solutions, and we can choose the smallest one).

\beT With   claims satisfying $\fz m_1=1$ (for example exponential), the readiness function is:
\bea f(\q) =\frac{(\Tl+\Tq)^2-\frac{K \Tl
   \Tq}{m_1}-\frac{\Tl}{m_1}}{\Tl
   \Tq}=\frac{\Tq}{\Tl}+\frac{2
   m_1-K}{m_1}+\frac{\Tl
   m_1-1}{m_1 \T q}.\eea

  Three possibilities arise for subsidiaries.  When $K=0,$ they are:
\BEN \im The  loading condition $\r <1$ is not satisfied, resulting in rejection.
\im The  loading condition $\r<1$ is  satisfied, and $k \leq \fr{(\l + \q) ^2 - \l c \mu}{\l \q} $ resulting in  acceptance.
\im The  loading condition $\r <1$ is  satisfied, and $k > \fr{(\l + \q) ^2 - \l c \mu}{\l \q} ,$ resulting in bailouts up to the level $\kil(\p_i),$ where $\p_i=\d(k) -\q,$ and $\d(k)$ is the solution of $k = \fr{(\l + \d) ^2 - \l c \mu}{\l \d} $.
\EEN
\eeT

  {\bf Proof}:  The derivative of the right hand side of \eqref{e:exef}  is proportional
to $\d^2 + \l(\mu c-\l)$, $f(\Tq,K)$  is increasing whenever the loading condition $\r <1$ is satisfied, and  therefore  an exponential subsidiary not satisfying the efficiency condition \eqref{e:exef} can always be rendered efficient by an appropriate killing, increasing its discount rate.

\section{Central branches as Sparre-Andersen processes\label{s:CB}}

Consider a  central branch with SA subsidiaries that
infuses capital into  subsidiary $i$  every time its surplus
level drops below $0$ for the $j$-th time. For simplicity, we assume from now on that the subsidiary is reset to zero. Given the challenges presented by   multidimensional CB networks, we
focus  now on   exponential or  \PH  subsidiary claims, and the latter are general enough for our purposes.

\begin{remark}    Taking $S_i(t)$ to be  spectrally positive Levy  perturbed compound  processes {\bf SPLPCP} \cite{frostig2008ruin,li2009renewal,zhang2013sparre} (for example with a Brownian motion perturbation)   poses
 often no problem.

In this case,   the  subsidiaries may be kept nonnegative  using  {\it minimal Skorohod regulation}. Then:
\begin{eqnarray}\label{dU0L}
&& X_0(t)= u_0 + c_0 t -S_0, \\&& S_0 =  \sum_{i=1}^I k_i I_i(t),  
\quad I_i(t) = -\inf_{s\leq t}\{X_i(s),0\}, \nonumber \end{eqnarray}
where the {\it regulator process} $I_i$ is the minimal
process whose addition to $X_i$  ensures that the sum is non-negative.
\end{remark}

With $(\vec \b,B)$ subsidiary claims,  the
  bailout (time, size) pairs $(\T A_k,\T C_k)_k$ are IID random variables
with joint distributions of the special form
\begin{eqnarray*}\label{2+}
P(A_k\in d t, C_k\in d x) = \vec \alpha(t) \, e^{ B_k x}\, \bff  b_k\,\, d  x\,d  t,\quad t,x\in\mathbb  R_+,
\end{eqnarray*}where $B_k:=k^{-1} B,$ and $\vec \alpha(t)=(\a_1(t), ...,\a_i(t), ...)$ contains the densities of the ladder time joint with ruin
in phase $i.$

The CB  is itself  a Sparre-Andersen process with
 phase-type claims,  exhibiting however a {\it non-standard} dependence \eqref{2I+}: the arrival times of this process are the first-passage times for the subsidiaries and the jump sizes are given by the bail-out amounts (undershoots).

 \begin{remark}  The classic Sparre-Andersen  model with  jumps of  phase-type $ (\vec \b,B)$  and independent  inter-arrivals with density $a(t)$ is a particular case of \eqref{2I+}, obtained  by taking densities of the form $$\vec \alpha(t)=a(t) \vec \b$$

 The essential difference in our case
  from the classic independent Sparre-Andersen process is that here the initial phase of a service(claim) period   is decided at the  bottom of the up-jump representing its inter-arrival  time,  and decides therefore also the size of the jump. \end{remark}

\begin{remark}\label{r:I>1} Consider a  CB network with several independent subsidiaries starting all at $u_i=0, , i=1,...,I$,  having claims of phase-type $\vec \b^{(i)}, B^{(i)}, i=1,...,I$,  and
let $\r^{(i)}(t), \bar R^{(i)}(t)$ denote the respective down ladder densities and survival functions.

By conditioning, we find that
 the density of the the first bailout of the CB  is
\begin{eqnarray}\label{2I+}
P(\T A_k\in d t, C_k\in d x) = \sum_{i=1}^I  \le( \prod_{j\neq i} \bar R^{(j)}(t)\ri)\, \vec \alpha^{(i)}(t) \, e^{ B^{(i)}_{k_i} x}\, \bff  b_{k_i}\,\, d  x\,d  t,\quad t,x\in\mathbb  R_+
\end{eqnarray}where $\vec \alpha^{(i)}(t)$ have Laplace transforms satisfying Kendall equations.

In the case of one subsidiary, $I=1$ the above equation gives the density for all bail-outs. When $I > 1$ this is no longer the case since  only one of the risk processes is controlled  at the time of  the first bailout.

The Sparre-Andersen process representing the CB  is considerably harder to analyze when $I > 1$.
In the sequel, we analyze the CB process in a stylized case of hierarchical networks for $I > 1$ and, respectively, we provide approximations of many quantities of interest based on the CB process in the case $I = 1$.

  \end{remark}

\section{Linear networks of deterministic CB: reduction to one dimension\label{s:red}}
It turns out that as long as the CB in isolation is a {\bf deterministic drift} with parameters $u_0,c_0$,  the {\bf ruin probability} of the CB  equals that of a subsidiary with  modified parameters
$$u_1'=u_1 +u_0/k,\quad c_1'=c_1+ c_0/k,$$
where  $k=k_1$, independent of the reset  policy!

The result is the same as if the CB transfers everything at the time $0_+$. This is  also the case with several other  problems involving a drift CB with no extra liabilities, which ends up  liquidated totally.

 {\bf The  pooled assets auxiliary process.\label{s:red}}
 Introduce
  the  pooled assets  combining $X_0$ and the reflected processes $\T X_i(t)=X_i(t)+I_i(t)$
  in such a way that the transfers and regulation  cancel out:
\beq X(t)&=& X_0(t)+ \sum_{i=1}^I k_i \le(X_i(t)+I_i(t)\ri)\nonumber \\&=&u_0 + c_0 t + \sum_{i=1}^I k_i X_i(t) \nonumber \\&=& u+   c t - \sum_{i=1}^I k_i S^{(i)}(t)   \nonumber \label{wres} \eeq
where we put $$u= {u_0} +  \sum_{i=1}^I k_i u_i, \quad c=  c_0 + \sum_{i=1}^I k_i c_i.$$

\begin{remark} Note that with spectrally negative Levy subsidiaries $X(t)$ is also a {\bf spectrally negative Levy process}, while $X_0(t)$ is  a complicated superposition of SA processes.

However, when $I=1,$  the ruin time of $X_0$ and $U$ coincide. Furthermore,
the  pooled reserves from the point of view  of the subsidiary
$$\fr{X(t)}{k_1}=u_1 + \fr{u_0}k + t(c_1 + \fr{c_0}k) -S^{(1)}(t)$$
 has the same law  as the subsidiary with  {\it combined initial value and income rate}!

 \end{remark}  %

\begin{remark} With several subsidiaries, \beq  \t=\inf\{t: X(t) <0\}.
\label{e:parg}
\eeq
represents the ruin time if the subsidiaries may start helping each other at no cost, once the CB is ruined. \end{remark}
These remarks yield the following:

\begin{theorem} \label{easy} Let $X_0$ be a CB with deterministic drift and {\bf arbitrary structure  subsidiaries}.

 i) Assume  $I=1$ and put $k=k_1$. Then, the ruin time $\tau_0$  of the MRN equals a.s. and in distribution the time $\t$ of the pooled process defined in
\eqref{e:parg},  and equals furthermore  the ruin time of
the subsidiary with modified initial reserve $u_0/k + u_1$  and premium
rate $c_0/k + c_1$,    $F_i(x)$,
namely \beq\label{psi_2st} \Rui(\bff u, \bff c,t) = \Rui_1(u_0/k +
u_1,c_0/k + c_1,t), \eeq
 independently of the reset  policy!

ii)
For  $I>1$, the time $\t$    is an {\bf upper bound} for the ruin time $\t_0$ of the CB:
$$ \t_0  \leq  \t.$$

iii)  Statements  similar  to A) hold for any Gerber-Shiu objective, with or without dividends
to the subsidiary, as long as the CB is a deterministic drift.

\end{theorem}

  {\bf  Optimal allocation  of total reserves $u_+=u_0+ u_1$ and premium rate $c_+=c_0+ c_1$}.
 Note that:
\begin{enumerate}
\item When $k \geq 1$,  $u_0=c_0=0, u_1=x, c_1=c$ achieve the {\bf  minimal ruin probability}.

 \im  For $k=1,$ the ruin probability is independent of the amount  {  $u_1 \in [0,u_+]$, as well as of the amount of premium rate $c_1 \in [0,c_+]$}.
 \end{enumerate}

  This optimization result fits the intuitively clear fact that with one subsidiairy and no
expenses, it is optimal to take advantage of the first transfer without cost to transfer  everything to the subsidiairy.

The next corollary extends the previous result to a linear chain of CB, which can represent a chain of reinsurers, see \cite{gerber1984}.
\begin{corollary} \label{c:easy} Let $X_0,...,X_{I-1}$ denote a {\bf linear chain} of CB's with deterministic drift.  Assume $X_i, i=0, ...,I-1$  must pay proportional costs $k_i$ for bailing out $X_{i+1}$.

  Then, the probability of ruin   of the MRN satisfies
   \begin{eqnarray*}\label{psi_h} \Rui(\bff u, \bff c,t) = \Rui_I(\fr{u_0}{k_0...k_{I-1}} + \fr{u_1}{k_1...k_{I-1}}+ ...
u_I,\fr{c_0}{k_0...k_{I-1}} + \fr{c_1}{k_1...k_{I-1}}+ ...
c_I,t), \end{eqnarray*}
where  $\Rui_I$ is the ruin probability of the   {\bf last subsidiary}
in the chain.
 \end{corollary}

\begin{remark} One interesting feature of this result  is that it does not require
any assumption on the probabilistic structure of the
subsidiary risk process.

 Another interesting feature is that similar reductions hold for other problems, as long as $I=1$  and  the main branch is a deterministic drift (in the absence of subsidiaries).  For example, one may add subsidiary dividends, ruin observed only
 at Poissonian
 times, Parisian ruin, etc. \end{remark}
\begin{Exa}Consider a CB with a {\bf spectrally two sided}  Levy subsidiary  $X(t)$, and {\bf ruin time $T_0$ observed only
 at Poissonian
 times $T_i, i=1, ...$ of frequency $\n$}. Recall \cite[(2)]{albrecher2015strikingly} that the survival probability in this case  is given by
 \begin{eqnarray*}\f_\n(u)=E \f(u+U), \end{eqnarray*} where $U$ is  the {\it up} factor of the Wiener-Hopf decomposition
 \begin{eqnarray*}U=\sup_{0\leq s \leq \kil_\n} X(t)=U(\kil_\n)- \inf_{0\leq s \leq \kil_\n} X(t)\end{eqnarray*}

 Supposing furthermore that the subsidiary  owes an {\bf exponential ruin severity penalty}, one may use
 \cite[(6)]{albrecher2015strikingly} that the survival probability in this case  is given by
 \begin{eqnarray*}E_{u,\n} \le(e^{- \q  T_0 + \th U(T_0)} ; T_0 < \infty \ri)= E\le[e^{- \q  T^U} E_{u +U} \le(e^{- \q  \t_0 + \th U(\t_0)} ; \t_0 < \infty \ri)\ri]  Ee^{- \q  T^D+ \th D}, \end{eqnarray*} where $D=\inf_{0\leq s \leq \kil_\n} X(t)$  is  the {\it down} factor of the Wiener-Hopf decomposition, and $T^U, T^D$ are the first times occurrences  of $U,D$.

 Note that dividends over a barrier, and/or Parisian ruin may  be added as well, using \cite[Prop. 4]{albrecher2015strikingly}, and \cite[Prop. 4]{albrecher2015strikingly}.

  If furthermore the  Levy subsidiary is spectrally negative,
the joint Laplace transform of the Poissonian observed ruin time
   and of the ruin severity
 is more explicit \cite[(14)]{AIZ}
 \begin{eqnarray*}&&E_u [e^{-q T_0 + \th U(T_0)} ; T_0 < \I]=\\&&\le( Z(u,\th) - Z(u, \F_\n)
\fr{(\F_\n -\F) \k(\th)}{\n(\th -\F)} \ri)( I - \n^{-1} \k(\th))^{-1}.
\end{eqnarray*}
Finally, the ruin probability of the CB will be given by the same formulas, with modified initial capital $u$ and premium $c$.
\end{Exa}

 Further questions of
interest,  not considered here,  are:

  \beQ Determining the  {\bf optimal allocation}  of the initial {\bf total capital} $u_+=\sum_{i=0}^I u_i$ and   {\bf total income rate} $c_+=\sum_{i=0}^I c_i$, when $ I \geq 2$, for example  under red time objectives.
  \eeQ

 \beQ Another  important   question is how to {\bf dynamically manage the subsidiaries}, i.e. how to set the reset  levels $\chi_j^{(i)}$ optimally,  possibly adaptively, depending on the evolution of the information available.

\eeQ

\beQ Study  CB branches with an extra Levy spectrally negative perturbation: for example with  exponentially distributed claims, or with Brownian motion (BM). The problem is particularly important, since the BM could model an approximation of  other subsidiaries, and is considered in the next section. \eeQ

\section{SNMAP approximations of non-deterministic CB}
In the reminder of the paper we focus on the case of a non-deterministic CB with one subsidiary. We have seen in the previous sections that the risk process of the CB is itself a Sparre-Andersen process in which the inter-arrival times are the first passage times and respectively the bailout amounts for the subsidiary.
We provide here a methodology to arrive at an SNMAP approximation the CB risk process, with the goal to provide an input to the available scale methodology toolbox (\cite{KP,Iva,IP, Ivapackage}) which can be used to simultaneously solve a wide variety of insurance problems for this CB risk process. We exemplify in this section with the computation of the optimal dividend barrier for the CB  itself, and implicitly for the whole network.
\subsection{Two point Pad\'e  approximations for  the downward ladder time of the \CL process
 with  exponential claims \label{s:rho}}

In this subsection we consider the particular case when the subsidiary risk process follows a \CL process
 with  exponential claims and we are interested in its downward ladder time.

 Our approach is based on the idea that approximating the excursions of a process ensures approximating the process, and in particular various functionals of the process \cite{lambert2013scaling,yano2013functional}.

  With phase-type jumps, one would need to provide {\bf bivariate matrix-exponential  approximations} for the {\bf joint law of the downward ladder time and height} of a SA process.

In the  case of  {\bf exponential claims},  the density of the downward ladder time  may be expressed as a  hypergeometric function:
\beq  \r(t)=\T \r(c \mu t) c \mu, \quad \T  \r(t) = \r e^{-(1+ \r) t}  \; _0 F_1 (2,  {\r} t^2 ).  \label{e:sc}\eeq

However, what we need is a
 {\bf  phase-type}  approximation of this.  This topic has already
 been considered in \cite{abate1988approximations} (at order two), as one of many possible methods for approximating the M/M/1  busy period density.

 We recall now some basic facts on this case \cite{AM15}.

 \begin{Lem} a) With Poisson arrivals of rate $\l$
and exponential claims of rate
$\mu $, the Laplace transform of the downward ladder time density satisfies a quadratic  equation
$$\Hr_\q  = \H a(\q  + c\mu - c\mu \Hr_{\d}= \fr{\l}{\l+\q  + c\mu - c\mu \Hr_{\d}}= \fr{\r}{1+ \de + \r -  \Hr_{\d}}, \; \r:=\fr{\l}{ c \mu}, \; \de=\fr q{c \mu }$$ with solution
\beq \label{e:LTe} \Hr_{\d}=\H {\T \r}_{\de}=\frac{1}{2 }\le(\de+\r +1  -\sqrt{(\de+\r +1)^2-4   \r
   } \ri).\eeq

 b) The Laplace transform \eqref{e:LTe}  may be computed iteratively by the continued fraction expansion
\begin{eqnarray*} \H {\T \r}_{\de}=\fr{\r}{1+ \de + \r- \fr{\r}{1+ \de + \r- \fr{\r}{1+ \de + \r-s...}}}.\end{eqnarray*}
 The  convergents $\H {\T \r}_{\de}^{(n)}, n=1,2,3...$ with $s$ constant satisfy   $\H {\T \r}_{\de}^{(n)}=\fr{\r P_{n-1}(\delta +\rho +1)}{P_n(\delta +\rho +1)}$ \cite[(75)]{neuts1966queue} where $P_n(x)$ are Chebyshev polynomials \cite[(77)]{neuts1966queue}. When $s=0,$ the first three are
\begin{eqnarray*} \left\{\frac{\rho }{\delta +\rho +1},\frac{\rho
   (\delta +\rho +1)}{\delta ^2+2 \rho  \delta +2
   \delta +\rho ^2+\rho +1},\frac{\rho  \left(\delta
   ^2+2 \rho  \delta +2 \delta +\rho ^2+\rho
   +1\right)}{(\delta +\rho +1) \left(\delta ^2+2 \rho
    \delta +2 \delta +\rho ^2+1\right)}, ...\right\}\end{eqnarray*}
\end{Lem}

Decomposing in partial fraction the third convergent
 yields  an order three rational approximation
 of the Laplace transform of the ladder time:
  \begin{eqnarray*} &&\Hr_{\d} \approx \frac{\r}{\delta + \r+1}
  \fr{(\delta + \r+1) ^2-\rho }{  (\delta + \r+1) ^2-2 \rho }=\\&&\frac{\r}{2}\le(\fr 1 {\delta+ \r+1 }+ \frac{1/2 }{\delta + \r+1 +\sqrt{2 \rho }}+\frac{1/2 }{\delta + \r+1-\sqrt{2 \rho }}\ri).\end{eqnarray*} 
Inverting the Laplace transform yields a hyperexponential density approximation:
 \beq \label{e:ratap} &&\T  \r(t) \approx \r e^{-(1+ \r) t}(\frac{1}{2 }+ \fr 14 (e^{-\sqrt{2 \r} t}+e^{+\sqrt{2 \r} t})):= \sum_{i=0}^2 \a_i \l_i e^{-\l_i t} \nonumber \\&& =\frac{\r}{2 }e^{-(1+ \r) t}+ \frac{\r}{4 } e^{-(1+ \r + \sqrt{2\r}) t}+ \frac{\r}{4 }e^{-(1+ \r - \sqrt{2\r}) t}, \eeq
 where $\a_0=\frac{\r}{2 (1+ \r)}$ and $ \a_{1,2}=\frac{\r}{4 (1+ \r \pm \sqrt{2\r})}$
  are nonnegative  for any $\r \in [0, \I)$. Furthermore,     $\sum_i \a_i <  \fr 34 \r$ iff $\r < 1$, providing us thus with a valid approximations for any $\r$ in this range.


\begin{remark}  A further simplification  of the Laplace transform \eqref{e:LTe} may be obtained factoring $\de+\r +1$ and changing variables $ a=\r (1+\r+\de)^{-2}$:
 \beq &&\Hr_{\d}=\H {\T \r}_{\de}=\frac{\de+\r +1}{2 }\le(1  -\sqrt{1-4   \fr{\r}{(\de+\r +1)^2}} \ri)\nonumber \\&&=\frac{\r }{\de+\r +1} \fr{1-\sqrt{1- 4 a}}{2 a}. \label{e:LT}  \eeq

The second factor put thus in evidence is
 the
   generating function of the famous Catalan numbers
   \begin{eqnarray*} \fr{1-\sqrt{1- 4 a}}{2 a}= \sum_{k=0}^\infty \fr{\binom{2 k}{k}}{k+1} a^k =1 + a + 2 a^2 + 5 a^3 + 14 a^4 +...\end{eqnarray*} and a  continued fraction (cf) representation
    \beq \label{e:cf} \fr{1-\sqrt{1- 4 a}}{2 a}=\fr{1}{ 1 -\fr{a}{ 1 -\fr{a}{ 1 +...}}}\eeq
   may be  found  for example in \cite[(7.7.5)]{cuyt2008handbook}.\fn[4]{
As well known,  the Pad\'e approximations obtained by truncating continued fractions have  good properties, like larger domains of convergence than the corresponding power series.}
The lowest order    approximations  \eqref{e:cf}  are
 $$\fr1{1-a}, \fr{1-a}{ 1-2{a}}, \fr{1-2 a}{ 1-3{a} + a^2}, \frac{1-3 a+a^2}{1-4 a+3 a^2},\frac{1-4 a+3 a^2}{1-5 a+6 a^2-a^3},...
 $$
\end{remark}

\begin{Lem}  The rational convergents
 $R_n$  of the continued fraction    \eqref{e:cf} increase
  towards $\fr{1-\sqrt{1- 4 a}}{2 a}, \forall a \in (0,1/4)$.
  \end{Lem}
  {\bf Proof:} This is immediate by  the positivity of $a$.

\begin{remark} Alternatively, we may use   two point Pad\'e approximations  which
 ensure also the equality of the derivatives around $0$.  \cite[Sec 3]{abate1988approximations} provide an in-depth  numerical comparison of several hyper-exponential approximations of order two, and find that fitting the derivatives yields excellent results around $0$, while fitting the moments is
 less satisfactory, since better results  may be obtained with asymptotic approximations. \end{remark}

 Let us invert now the lowest order  two-point Pad\'e  approximation of $\Hr(\d)$ which
 ensures also the condition $\T {\Hr}_{\d=0}=\r$, $\r \in (0,1]$: 
\beq  \label{e:tpp} &&\frac{\rho  \left(\delta ^2+2 \rho  \delta +2 \delta +\rho +1\right)}{(\delta +\rho +1) \left(\delta ^2+2
   \rho  \delta +2 \delta +1\right)}=\fr {\r}{\rho +2}\le(\frac{\rho +1 }{ \delta +\rho +1}+\frac{\delta   + \rho +1 }{ \delta
   ^2+2   \delta(\rho +1) +1}\ri) \\&&=\fr {\r}{\rho +2}\le(\frac{\rho +1 }{ \delta +\rho +1}+\fr 1{\lambda _1-\lambda _2}\le(\frac{{\lambda _1-(\rho +1)}}{ \delta
   +\l_1}+\frac{{\rho +1-\lambda
   _2}}{ \delta+
   \l_2}\ri)\ri), \l_{1,2}=\rho +1 \pm \sqrt{\r(\r+2)}\nonumber     \eeq

This yields a  density approximation:
   \beq \label{e:rat2p} &&\T  \r(t) \approx \sum_{i=0}^2 \a_i \l_i e^{-\l_i t}, \eeq
 where $\a_0=\frac{\r}{ \r+2}, \a_1=\frac{\r}{ \r+2} \frac{{1-(\rho +1)/\l_1}}{\lambda _1-\lambda _2}=\frac{\rho ^{1/2} \left(\sqrt{\rho  (\rho +2)}+1\right)}
 {2 (\rho
   +2)^{3/2} \left(\rho +1+ \sqrt{\rho  (\rho +2)}\right)},$

   $\a_2=\frac{\r}{ \r+2} \frac{{(\rho +1)/\lambda
   _2-1}}{\lambda _1-\lambda _2}=\frac{ \rho ^{1/2}\left(  \sqrt{\rho
   (\rho +2)}-1\right)}{2 (  \rho +2)^{3/2} \left(\rho +1-\sqrt{\rho  (\rho +2)}\right)}$.  Note that $\a_i$
  are nonnegative for any $\r \in [0, \I)$, and $\sum_i \a_i =\r$.



 \begin{figure}[h!]
\begin{center}
\includegraphics[width=0.45\textwidth]{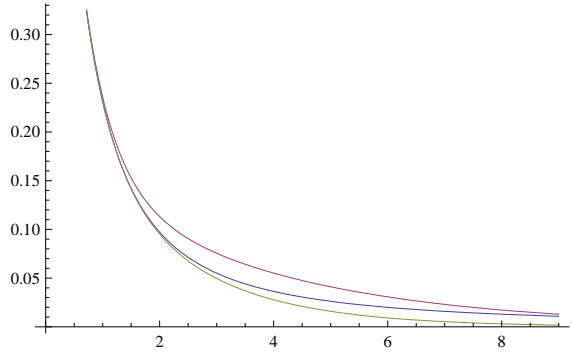}
\end{center}
\caption{Two order 3 matrix exponential approximations
\eqref{e:ratap}, \eqref{e:rat2p} of $ \r(t)=\rui(t;0) $  with $\l =\mu=1, \r=\fr 78=c^{-1}$. The Bessel law (blue) is in between
  the two point Pad\'e approximation (in red), which is tight both at $0$ and $\I,$  and the classic
  continued fraction Pad\'e approximation (in yellow), which is better near $0$}
\label{f:Pade}
\end{figure}

\begin{figure}[h!]
\begin{center}
\includegraphics[width=0.45\textwidth]{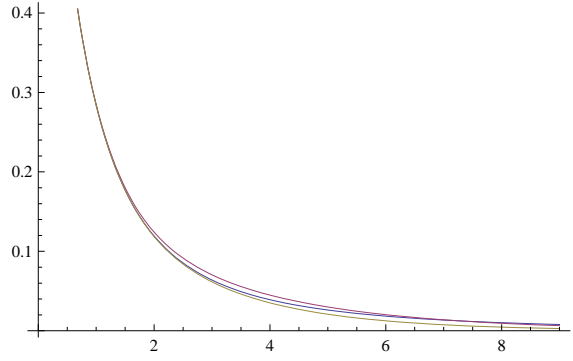}
\end{center}
\caption{Two order 3 matrix exponential approximations
\eqref{e:ratap}, \eqref{e:rat2p} of $ \rui(t;0) $  with $\l =\mu=1, \r=\fr 12=c^{-1}$. The fit is better as $\r$ gets further from $1$, and could be improved by using a combination of the two approximations, of order $5$.}
\label{f:Pade}
\end{figure}


%
%
%
%
%
%

\subsection{Pad\'e  based SNMAP
approximations for $X(t)$  when subsidiary claims are exponential
\label{s:PDV}}
{


Having established in the previous subsection convenient approximations for the subsidiary risk process, we now approximate the CB risk process.

Note that even though finite time ruin probabilities have an  explicit Bessel density
with exponential claims, we will  replace them by matrix exponential approximations, since this allows solving network problems
by the SNMAP methodology.


  After applying the order three  approximation \eqref{e:rat2p} to the  subsidiary's ladder time, the central branch becomes  a MAP with three states, with transition
   rates $Q_{ij}= \l_i \a_j$ accompanied by exponential jumps of rate $\mu/k$ translated by $K$ (we could include here {\bf phase-type jumps to the CB}, and {\bf fixed costs}, since these pose no  problem to the MAP methodology).

   When $K=0,$ $\H f(s)=\frac{\mu/k  }{s+\mu/k },$ and we find from the general formula  that  the  symbol of the approximated CB
   is \beq \label{e:MAPP} &&\K(s)=diag(c_0 s - \l_i)+ \bff \l \vec \a \H f(s)\\ &&= \left(
\begin{array}{lll}
 c_0 s-\l_0 +\l_0 \a_0 \frac{\mu/k  }{s+\mu/k } & \l_0
   \a_1 \frac{\mu/k  }{s+\mu/k } & \l_0 \a_2 \frac{\mu/k  }{s+\mu/k }\\
 \l_1 \a_0 \frac{\mu/k  }{s+\mu/k } & c_0 s-\l_1 +\l_1 \a_1
    \frac{\mu/k  }{s+\mu/k } & \l_1 \a_2 \frac{\mu/k  }{s+\mu/k } \\
 \l_2 \a_0 \frac{\mu/k  }{s+\mu/k } & \l_2 \a_1 \frac{\mu/k  }{s+\mu/k }
   & c_0 s-\l_2+\l_2 \a_2 \frac{\mu/k  }{s+\mu/k }
\end{array}
\right). \nonumber \eeq

To apply the scale based MAP methodology, it is convenient
 to transform this MAP with exponential jumps into a continuous MMBM:

\beq \label{e:MMBM} \T {\K}(s)=\left(
\begin{array}{llllll}
 c_0 s-\l_0 & 0 & 0 & \l_0 \a_0 & \l_0 \a_1 & \l_0 \a_2 \\
 0 & c_0 s-\l_1 & 0 & \l_1 \a_0 & \l_1 \a_1 & \l_1 \a_2 \\
 0 & 0 &c_0 s-\l_2 & \l_2 \a_0 & \l_2 \a_1 & \l_2 \a_2 \\
 \mu/k  & 0 & 0 & -s-\mu/k  & 0 & 0 \\
 0 & \mu/k  & 0 & 0 & -s-\mu/k  & 0 \\
 0 & 0 & \mu/k  & 0 & 0 & -s-\mu/k \end{array}
\right) \eeq

After obtaining an SNMAP approximation for $X(t)$, we may solve approximatively various problems related to this process, using the package \cite{Ivapackage}.


We illustrate this by finding in Figure \ref{f:Div} the  optimal dividend barrier. It is beyond the scope of this article to investigate further this problem in particular (which deserves a separate treatment), as our main goal was not the solution of the  the problem per se but providing an adequate input to the toolbox that solves this type of problems.

\begin{figure}[h!]
\begin{center}
\includegraphics[width=0.5\textwidth]{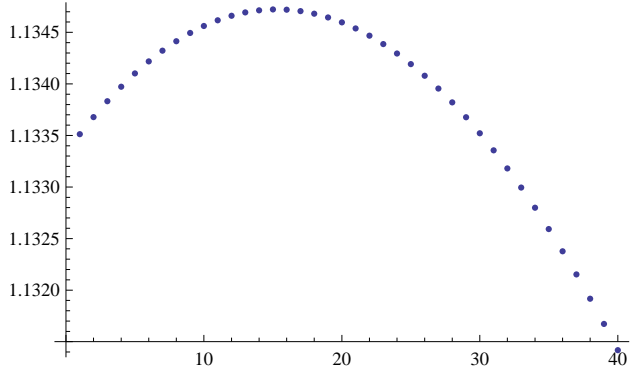}
\end{center}
\caption{Dividends as function of barrier, starting from $u_0=u_1=0$, with $\l=\mu=k=1,$ $c_1=4, c_0=24$, obtained by the scale matrix methodology \cite{Iva}. The maximizing barrier is $\approx .17$}
\label{f:Div}
\end{figure}

\section{Spectrally-negative Markov additive processes  (SNMAP)- a new framework for risk modelling}\label{s:MAP}
Having demonstrated the power of the approximation based methodology, we dedicate the remainder of the paper to set bases of future research directions that build on the extraordinary potential in insurance of the class of spectrally-negative MAPs  (SNMAP),
that is, MAPs which have negative jumps only.

 Informally, a MAP is a bivariate process $(X(t),J(t))$,
where $J(t)$ is a Markov chain (MC) representing an exogenous background process, and $X(t)$
is a so-called additive component modulated by $J(t)$ (nonetheless we often say MAP meaning $X(t)$).

A MAP is a generalization of a L\'evy process in the sense that $X(t)$ has stationary and independent increments conditioned on the state of
the modulating process $J(t)$.
\begin{Def}\label{def:MAP}
A bivariate process $(X(t),J(t))$ is called MAP if, given $\{J_{T}=i\}$, the shifted process
 $(X_{T+t}-X_{T};J_{T+t})$ is independent from $(X(t),J(t)),0\leq t\leq T$ and has the same law as $(X(t)-X(0);J(t))$ given $\{J(0)=i\}$ for all $i$ and $T>0$.

  A spectrally-negative MAP is a MAP whose additive component can have only negative jumps.
\end{Def}

\begin{remark} Furthermore, one can replace the deterministic $T$ in the above definition by a stopping time. The resulting property is called the strong Markov property for MAPs.\end{remark}
\noindent It is common to assume that $J(t)$ is an irreducible MC with a finite state space $
\mJ=\{1,\ldots,n\}$, which we do throughout this work.
It can be shown that $X(t)$ evolves as some L\'evy process $X^i_{t}$ while $J(t)=i$, and in addition a transition of $J(t)$ from $i$ to $j$ may trigger a jump of $X(t)$ distributed as $U_{ij}$,
where $J(t)$ and all the components in the construction are assumed to be independent. This construction presents an alternative often-used definition of a MAP. See~\cite[Ch.~XI]{APQ} for further information on MAPs.\\

\begin{remark} In the presence of a "modulating environment" $J$, one  key idea   is to consider
 matrices of expectations and probabilities, conditioned on some phase $i$ at start and joint with some phase $j$ at crossing the level of interest. These will be denoted by $\E=\E[... ;J(t)], \P=\P[... ;J(t)]$. \end{remark}
 {\bf The Laplace exponent}. The law of a spectrally-negative MAP is characterized by a certain matrix-valued function $\K(\theta)$,
 defined by
 \[\E[e^{\theta X(t)};J(t)]=e^{\K(\theta)t},\]
where the $(i,j)$-th element of the matrix on the left is given by $E(e^{\theta X(t)};J(t)=j|J(0)=i)$
\fn[4]{This is the analogue of the "Laplace exponent/symbol/cumulant generating function" of a spectrally-negative L\'evy process.}.
The $n\times n$ matrix $\K(\theta)$ is given by
\bea \K_{ij}(\theta)=\begin{cases}
  \k_i(\theta)+q_{ii},&\text{if }i=j,\\
  q_{ij}\H U_{ij}(\theta),&\text{if }i\neq j,
\end{cases} \eea
where $q_{ij}$ are the elements of the transition rate matrix $Q$ of $J(t)$, $\k_i(\theta)=\log \E e^{\theta X_i(1)}$ is the Laplace exponent of the L\'evy process $X^i_{t}$, and
$\H U_{ij}(\theta)=\E e^{\theta U_{ij}}$. In matrix form
\begin{eqnarray*}\K(\th)=Diag(\k_i(\th))+ Q \o \H U(\th) \end{eqnarray*}where $\o$ denotes Hadamard product.

\begin{remark} The Levy processes may have  killing parameters $\q_i=-\k_i(0) \geq 0$. When $\k_i(0)=0$,  it follows that $\K(0)=Q$  is the transition rate matrix of $J$. \end{remark}
\begin{remark} An alternative  notation in use when a fixed nonzero killing $\q$ is involved is to   exclude it  from the symbol matrix $\K (\th)$, and  write $\K_{\q}(\th)=\K (\th)-\q I$ when killing  is present.
We will follow, for the sake of readability, the convention of  including a fixed  killing argument $\q$ in the symbol, scale, etc
   and just write $ \K,W,...$ instead of  $ \K_{\q}, W_{\q},...$.

 If an additional secondary killing $\n$ appears  in certain circumstances,
 its presence must be indicated. The notations  $ \K_{\p}, W_{\p},...$   corresponds then to
 $ \K_{\p +\q}, W_{\p +\q},...$ in the complete notation.

\end{remark}

\begin{Exa}
A Markov-modulated Cram\'er-Lundberg process with tax  is retrieved from a spectrally-negative MAP $(X(t),J(t))$ by putting $U_{ij}=0$ and
$X^i_{t}=c_i t-\sum_{k=1}^{N_i(t)} C_k$, where $N_i(t)$ is a Poisson process of intensity $\lambda_i$ and $C_k$ are iid positive random variables.
Hence $G_{ij}(\theta)=1$ and $\k_i(\theta)=c_i\theta-\lambda_i(1-\E e^{-\theta C_1})$.
\end{Exa}

\begin{remark} With phase-type jumps, it may be more convenient to work with the symbol  of the "embedded process", in which  jumps have been replaced by downward drifts \cite{asmussen1995stationary,AA}.
\end{remark}

 In applications, we are interested  in "regulated" versions  of  $X(t)$:
 \begin{eqnarray*}&& X_{[0}(t)= X(t) + L_0(t), \; \;  X^{b]}(t)= X(t)  -  R_b(t),\\&& X^\g(t)=X_{[0}^{b],\g}(t)= X(t) + L(t) - \g R(t), \; \g \leq 1, \end{eqnarray*} where $L_0(t)=-(\und X(t) \wedge 0), \; \und X(t)=\inf_{0 \leq s\leq t} X(t)$ is the minimal  regulator constraining $X(t)$ to be nonnegative,
 $$R_b(t)= \le(\bar X(t)-b\ri)_+, \; \bar X(t):= \sup_{0 \leq s\leq t} X(t) $$
  is the minimal  regulator constraining $X(t)$  to be smaller than $b$.   The process    $X_{[0}^{b],\g}(t)$ is said to be  reflected at $0$ and  refracted  with coefficient $\g$ \cite[(1)]{AIpower}.   The  minimal  regulators $L(t),R(t)$ have points of increase contained in $\{t\geq 0:Y(t)=0\}$ and $\{t\geq 0:Y(t)=\bar Y(t) \vee b\}$ respectively, and their explicit formulas are more complicated  in this case; for a recursive construction, see for example \cite[pg 165]{APP}.  The process  $X_{[0}^{b],1}(t)$
  will be denoted  by $\T X^b(t)$.

  \beR
  The process $\g R(t)$ intervenes here as cumulative dividends paid to some beneficiary, and using $\g <1$  allows
  $X^\g(t)$ to continue above the upper "threshold" $b$ (with a modified drift). \eeR

\qq {\bf Smooth exit upwards and the associated generator}.
One key idea in the  study of SNMAPs is considering  the phase process "while progressing upwards" $J_{\tau_x^+}$, which is itself a MC as a function of $x\geq 0$, by Definition~\ref{def:MAP} with $T=\tau_x^+$.
The transition rate matrix  $G$ is   a (sub)generator,   the matrix
\beq \la{onesided} \P[J_{\tau_x^+}]=e^{G x}\eeq
contains the  probabilities of hitting $x$ in various phases,
and $$e^{G x} \bff 1= P[\tau_x^+ <\I]$$
contains the  probabilities of ever hitting $x$, starting from all initial states.

The matrix $G$ is a right solution of the equation $\K(-
G)=0$ \cite{breuerPH}.
 \begin{remark} In the Levy literature, $G$ is a scalar denoted by $-\F$, which solves the Cram\'er-Lundberg equation $\k(\F)=0, \F \geq 0$.\end{remark}
 If the MC $J_{\tau_x^+}$ is transient, then
 its life time coincides with the overall supremum of $X(t)$, and has thus a PH distribution characterized by the matrix $G$.

 {\bf The   two-sided smooth exit problem and the scale function}.
Cf. \cite{KP,IP}, the solution of the smooth two-sided SNMAP exit problem has a multiplicative form (well-known in the case of Levy  processes):

\begin{theorem}\label{thm:main_scale}
There exists a unique
continuous function
$W:[0,\infty)\rightarrow\R^{n\times n}$
such that \[\int_0^\infty e^{-\theta x}W(x)\dd x=\K(\theta)^{-1}\]
for all sufficiently large $\theta>0$, $W(x)$ is invertible for all $x>0$, and
\beq\label{s} \P_x[\tau_b^+<\t ^-_0;J(\tau_b^+)]=W(x)W(b)^{-1}\text{ for all }b\geq 0,  x \in [0,b].\eeq

The decomposition
\begin{eqnarray*}W(x)=e^{-G x} H(x),\end{eqnarray*}where $H(x)$ is the matrix of expected occupation times at $0$ up to the first passage over $x$, provides a probabilistic interpretation of the scale matrix.
\end{theorem}

\beR Informally, $W(x)$ is the solution of $\bc G W(x)=0, & x \geq 0\\ W(x)=0,& x <0\ec$, where $G$ is the Markovian generator of the process $X(t)$. \eeR

\begin{remark} A matrix expansion for   the two-sided exit  by jump probability with phase-type claims and inter-arrivals was already obtained in \cite[(13)-(14)]{APU} via an elementary ODE approach, later extended in  \cite[(26)-(28)]{jiang2008perpetual}. However, the intuition and rigorous proof that everything could be expressed in terms of one { scale matrix} were  furnished only later, in the seminal  papers \cite{KP,IP}, respectively. \end{remark} \bigskip
\begin{Exa}\label{ex:SA}
A Sparre Andersen renewal risk process with inter-arrival times of \PH $(\vec{ \a}, A)$,  is an example of a spectrally-negative MAP $(X(t),J(t))$.
Here the MC $J(t)$ lives on $n$ states, and  $J(0)$ is distributed according to $\vec{ \a}$. For $t >0,$ $J(t)$
makes a jump from $i$ to $j$ without causing a jump of $X(t)$ with rate $a_{ij}$; it makes a jump from $i$ to $j$ and triggers a jump $-C_k$ with rate $a_i\a_j$.
Hence for $i\neq j$ it holds that $F_{ij}(\theta)=a_{ij}+a_i\a_j\E e^{-\theta C_1}$, that is $q_{ij}=a_{ij}+a_i\a_j$ and $U_{ij}$ is an appropriate mixture of $0$ and $-C_1$.
Then $F_{ii}(\theta)=a_{ii}+a_i\a_i\E e^{-\theta C_1}+c\theta$, because $q_{ii}=-\sum_{j\neq i}q_{ij}=a_{ii}+a_i\a_i$ and $\k_i(\theta)=c\theta-a_i\a_i(1-\E e^{-\theta C_1})$ which corresponds to a compound
Poisson process with intensity $a_i\a_i$, jumps distributed as $-C_1$ and drift $c$. In matrix notation we have
\[\K(\theta)=A+\bff a\vec{  \a}\E e^{-\theta C_1}+c\theta\matI_n,\]
where $\matI_n$ is an $n\times n$ identity matrix.

\end{Exa}

 {\bf The   second scale function \cite{APP15,IP}} is defined by
 \begin{eqnarray*}Z(x,\th)= e^{\th x} \le( I- \int_0^x e^{- \th y} W(y) dy \; \K(\th) \ri). \end{eqnarray*}

 Alternatively,  defining
 the Dickson-Hipp transform
\begin{eqnarray*}\H W_x(\th)= \int_0^\I e^{- \th y} W(x+y) dy =e^{\th x} \le( \K(\th)^{-1}- \int_0^x e^{- \th y} W(y) dy \;  \ri), \end{eqnarray*}{ where the  equality holds for $\Re(\th)$ large enough},
it holds that \begin{eqnarray}Z(x,\th)=\H W_x(\th)  \K(\th) \Eq \H W_x(\th)=  Z(x,\th) \K(\th)^{-1} , \end{eqnarray} i.e the second scale function $Z(x,\th)$ coincides up to a constant matrix  with the Laplace transform of the shifted scale function (the "normalization" ensures that  $Z(0,\th)=I$).

Another way to characterize the second scale function is via its Laplace transform: \begin{eqnarray*}\H Z(s,\th)=(\th -s)^{-1}(\K(s)^{-1} -\K(\th)^{-1}) \K(\th)=(\th -s)^{-1}(\K(s)^{-1} \K(\th) -I).\end{eqnarray*}
In the Levy case, this becomes
\begin{eqnarray*}\H Z(s,\th)= \fr{\k(s)-\k(\th)}{(s-\th )\k(s)}.\end{eqnarray*}

{\bf The   second scale function with two-step killing} is \cite{AIrisk}, \cite[(12)]{AIZ}:
 \beq
 \la{e:s2}
 \P_x[\tau_b^+<T ^-_0;J(\tau_b^+)]=Z(x, \F_\p) Z(b, \F_\p)^{-1} \forall b\geq 0, x \in [0,b],
 \eeq
 in the sense that it replaces the first scale function in  the smooth two-sided exit problem  with extra killing below $0$.

\begin{remark}
The function $Z(x,\th)$ appeared first in 2008 in \cite{APP15}, in the context of finding smooth Gerber-Shiu functions associated to an {\bf exponential payoff} $e^{\th x}$.
Subsequently,  its key  role  in  many other first passage problems was revealed in \cite{IP} and  \cite{AIZ}.  \end{remark}

{\bf Killing when the draw-down exceeds a given value $a$}.
Consider    a process $X^{b]}(t)=X(t)-R_{b}(t)=X(t) -(\overline X(t) - b)_+$  regulated and started at  $b,$   resulting  in $$X^{b]}(t)= X(t) + b -\overline X(t)=b-Y^0(t),$$
where $$Y^x(t)=x + \overline{ X}(t) - X(t), \;  \;  Y^x(0)=x$$
 is called "drawdown/reflection from the running maximum" (starting from $x$) \cite[pg. 248]{K}.

The regulator $R_{b}(t)=(\overline X(t)-b)_+$ can also be seen as the total amount of dividends paid until ruin in a L\'evy model with the barrier dividend strategy,
where the initial capital and the barrier are both placed at the level~$b$.

Let us   kill now $X(t)$ (send it to some absorbing state) at the stopping time
\[t_a=\inf\{t\geq 0:Y_t> a\},\]
i.e.~at the first time when the {\bf height of an excursion from the maximum} exceeds~$a>0$, or, equivalently, at the first time when the regulated process $X^{b]}(t)$ started from $b>0$  drops by more than $b-a$.

Using the strong Markov property for MAPs, we see that also in the presence of killing at $t_a$, the environment phase observed while evolving upwards $J_{\tau_x^+},x\geq 0$ is still a MC, with some transition rate matrix $\L(a)$, so that
$\P[\tau_x^+<t_a;J_{\tau_x^+}]=e^{\L(a)x}.$

It was shown in~\cite[4]{IP}, generalizing  a well-known excursion theory relation for L\'evy processes, that for $a>0$ the right and left derivatives $W'(a_+)$ and $W'_-(a)$ exist and
\begin{align}\label{excscale}
 &\L(a)=-W'(a_+)W(a)^{-1},&\L(a-)=-W'_-(a)W(a)^{-1}, \; \; \P_x[\tau_b^+<t_a;J_{\tau_b^+}]=e^{\L(x+a)(b-x)}.
\end{align}

\begin{remark}
 For    SNMAP's, the life time of the transient MC characterized by $\L(a)$ is a matrix exponential law \eqref{e:divexp}.
 In the case of a L\'evy process, this reduces to an exponential random variable with rate $W'(a_+)/W(a)$ \eqref{e:divexp}.
    \end{remark}

\ssec{Resolvents}

The \emph{$q$-resolvent measure}, for any Borel set $B\in[a,b]$, may be expressed in terms of the scale function
\begin{align} \label{resolvent_density}
\begin{split}
	\E_x \Big( \int_0^{ \tau^-_a \wedge \tau_{b}^+ } e^{-qt} 1_{\left\{ X_t \in B  \right\}} \diff t\Big) &= \int_{a}^{b} 1_{\{y \in B\}}\Big[ \frac {W_q(x-a) W_q (b-y)} {W_q(b-a)} -W_q (x-y) \Big] \diff y \quad a < y < b;
	\end{split}
\end{align}
see Theorem 8.7 of \cite{K}.  { Letting $a \to -\I$ yields}
\begin{align}\label{resolvent}
\E_x\left(\int_0^{\tau_b^+}e^{-qt}1_{\{X_t\in A\}}dt\right)=\int_A\left\{e^{-{\Phi_q (b-x)}}W_q(b-y)-W_q(x-y)\right\}dy.
\end{align}

{\bf Resolvent of doubly reflected L\'{e}vy processes}  As
shown in (\cite[Thm. 1]{PDR}),  a version of the
$\q$-potential measure of a doubly reflected L\'{e}vy process with upper barrier $b$ is $\T U^b_\q(  x,d  y) = \int_0^{\I}e^{-qt}\P_x(\T X^b(t)\in d y)$
of $\T X^b$ is given by
\bea \T U^b_\q(  x,d  y) = \T u^q(x,b)\de_b(d  y) +
\T u^q(x,y)d  y\eea
 where $\de_b$ is the point-mass at $b$ and
\begin{equation}\label{eq:resolut}
 \begin{cases} \T u^q(x,y) = Z_\q(  x) Z'_\q(  b)^{-1}  W_\q'(b-y)
- W_\q(  x-y),  &x,y\in[0,b], y\neq b\\
\T
u^q(x,b)=Z_\q(  x) Z'_\q(  b)^{-1} W_\q(0),& y= b. \end{cases}
\end{equation}

The next section illustrates further the fact that the answers to a large variety  of first passage problems for SNMAP's may ergonomically be expressed  in terms of the matrices $\K,G,H,W$ and $Z$.

\section{A compendium of first passage formulas for SNMAP and SNL\'evy processes \label{s:MAPfor}}
In the research avenue we envision, diverse problems in insurance and in particular multidimensional problems which so far have defied exact analysis can be solved using approximations by SNMAP and "the eight pillars of one-dimensional first passage problems for SNMAP".
We  assemble here  these key formulas,  which have appeared in various prior works \cite{APP,AIZ,AIpower,albrecher2015strikingly} for the Levy case and \cite{IP,Iva} for the SNMAP case.

\subsection{Homogeneous problems}

\BEN

\im {\bf The De Finetti expected  discounted dividends over a barrier $b$},
satisfy  a relation  similar to \eqref{s}: \begin{eqnarray} \E_x\le[\int_{[0,\tau_{0}^-]}e^{-qs}d   R_b(s)\ri]=W_\q(  x)
W_{\q}^ {\prime}(b)^{-1}, \label{div} \end{eqnarray}
when $\t_{0}^-$ is the classic ruin time.

With {\bf Poissonian observed ruin}, the expected discounted dividends over a barrier $b$ become
  \cite[(27)]{AIZ} \begin{eqnarray} \label{div2} \E_x\le[\int_{[0,T_{0}^-]}e^{-qs}d   R_b(s)\ri]=Z_\q(  x, \F_{\q +\n})
Z_{\q}^ {\prime}(b, \F_{\q +\n})^{-1}, \forall x\in  [0, b]. \end{eqnarray}
 In the case of continuous observations $\p \to \I$,  this expression reduces to the previous result \eqref{div}.

The 
law of the total dividends  until ruin $R_b(\tau_{0}^-)$ is known as well; it is a matrix exponential distribution
generated by $\L(b)=-W'(b)W(b)^{-1}$ -- see \eqref{e:divexp}.

\im {\bf Expected  discounted dividends until the total bail-outs of a  reflected  process surpass an exponential variable $\kil_\th$}
satisfy (see \cite[(15)]{AIpower} for the particular case $b=x$)
 \begin{eqnarray}V= \E_x\le[\int_{[0,\I]}e^{-\q s} 1_{ [L_0(s) < \kil_\th]} d  R_b(s) \ri]=Z_\q(  x,\th)
Z_{\q}^ {\prime}(b,\th)^{-1}. \label{e:divskillbo}\end{eqnarray}

 When $\th=0$, this yields \cite[(4.3)]{APP}, and
when $\th \to \I$, this recovers \eqref{div}.

{}

\im {\bf Bail-outs of a  reflected  process, until the first dividend}.
Let $X_{[0}(t)$ denote a SNMAP process reflected at $0,$
 let $L_0(t)=-(0\wedge \underline{X}(t))$ denote its regulator at $0$, so that $X_{[0}(t)=X(t) +L_0(t) $, and let $\E^0_x$ denote expectation for the process   reflected at  $0$. Then, the joint law of the time until a process reflected at the infimum  climbs to an upper level $b$ and of the  bail-out (regulation) is:
\beq \label{refbailouts} \E^0_x [e^{-\q \t_b^+ - \th L_0(\t_b^+)}]=\bc Z_{\q}(x,\th) Z_{\q}(b,\th)^{-1} & \th <\I\\\P[ \t_b^+ < \t_0^-]=
W_{\q}(x) W_{\q}(b)^{-1}& \th =\I\ec,\eeq
where $Z_{\q}(x,\th)$ is the second scale function \cite[Thm 2]{IP}.

Furthermore, \cite[(3)]{AIpower} show in the Levy case that a power relation holds when replacing  $X_{[0}(t)$  by a refracted process $Y_0^{b,\g}(t) $ started at $b$. The proof uses the  probabilistic interpretation
\bea \E^0_x [e^{-\q \t_b^+ - \th L_0(\t_b^+)}]=P[\t_b^+ < \kil_\q
\wedge K_{\th}],\eea
where $K_{\th}$ is the first moment when the total bail-out exceeds
an independent exponential rv. $\kil_\th$.

 Finally,  \cite[(22)]{AIZ}  extend to the case when $\t_b^+$ is replaced by $T_b^+$.

{\iffalse \color{red} When $b \to\I$,  \eqref{refbailouts} yields ? \fi}

\EEN

\subsection{Non-homogeneous problems}

\BEN
\im {\bf The severity of ruin before seeing a barrier $b$ (two-sided exit)}.
Applying  the previous  result \eqref{refbailouts} and the two-sided exit formula \eqref{s},  one finds that \cite[Cor 3]{IP}:
\begin{eqnarray*}&& \E^0_x [e^{-\th L_0(\t_b^+)}; J(\t_b^+)]=Z(x,\th) Z(b,\th)^{-1}=
\\&&\P_x[\t_b^+ < \tau_0^- ; J(\t_b^+)] +
\E_x [e^{\th X(\tau_{0}^-)} ;  \t_0^- < \t_b^+ ] \; \E^0_0 [e^{-\th L_0(\t_b^+)}; J(\t_b^+)]\\&&= W(x)W(b)^{-1} + \E_x [e^{\th X(\tau_{0}^-)} ;  \tau_0^- <\t_b^+ ] Z(b,\th)^{-1}.\end{eqnarray*}

 We may solve now for
the joint Laplace transform of the first passage time of $0$, and the undershoot
\begin{equation} \label{sevruin}
\E_x\left(e^{-\q\t_0^- + \th X(\t_0^-)};\t_0^-<\t_b^+\right)=Z_\q(x,\th) - W_\q(x) W_\q(b)^{-1} {Z_\q(b,\th)} \end{equation}

In the SNMAP case, denote by $$H=\lim_{x\to \I} H(x)$$
the matrix of total expected occupation times at $0$.  Assume that either $Q \bff 1 \neq 0, $ or $\k'(0) \neq 0$.

Letting $b \to \I$ above, and using \cite[Cor 4]{IP} \beq\label{e:lWZ}  \lim_{b \to \I} W(b)^{-1} Z(b,\th)= (R+ \th I)^{-1}  {\K(\th)}:=\T \K(\th), \eeq   where $R=H^{-1} G H$ is a left solution of the equation $\K(-
R)=0$,  we find the {\bf  severity of ruin (non-smooth one-sided exit)}:
\begin{eqnarray*}&&\;  \forall x,\th \geq 0, \quad \E_x [e^{\th X(\tau_0^-)} ; \tau_0^- < \I]=Z(x,\th) - W(x)   (R+ \th I)^{-1}  {\K(\th)} \\&& = Z(x,\th) - W(x) H^{-1} (G+ \th I)^{-1} H {\K(\th)}. \end{eqnarray*}

\begin{remark} Note that if $\th$  is a zero of $det(\K(\th))$, then it is a zero of $det(\th I + G)$ and $det(\th I + R)$.  Thus, $\T \K(\th)$ should
be interpreted in a limiting sense. 
\end{remark}

In the Levy case, \eqref{e:lWZ} reduces to $\lim_{b \to \I} W(b)^{-1} Z(b,\th)= \fr {\k(\th)}{ \th-\F} $, yielding  the classic \cite[(7)]{AIZ}
\beq\label{e:old}
&&\E_x [e^{-\q \tau_0^- +\th X(\tau_0^-)} ]= Z_\q(x,\th) - W_\q(x)  \fr {\k(\th)-\q}{ \th-\F_\q}. \eeq

\beR When we  apply the above with $\th=\Phi_\p $, this becomes:
\begin{align}\label{sevruinb} &&
\E_x\left(e^{-\q\tau_0^- + \Phi_\p X_{\tau_0^-}};\tau_0^-<\infty\right)
=Z_\q(x,\Phi_\p ) - W_\q(x)  \fr {\p -\q}{ \Phi_\p -\Phi_\q}= \nonumber \\&&e^{\Phi_\p x}\left(1+(\q-\p)\int_0^xe^{-\Phi_\p z}W_\q(z)dz-\frac{\q-\p}{\Phi_\q  -\Phi_\p }e^{-\Phi_\p x}W_\q(x)\right).
\end{align}
\eeR
\im {\bf The severity of ruin at a Parisian ruin time, before seeing a barrier $b$}.  When $\t_0^-$ is replaced by a Parisian ruin time $\t_\p$, the severity of ruin satisfies in the Levy case \cite[(15)]{AIZ}, \cite[(1.13)]{BPPR}
\begin{equation}\label{extL}
\E_x \Big[ \mathrm{e}^{-\q \pt + \th X(\pt)} , \pt < \tau_b^+ \Big] =\left(Z_{\q}(x, \th) -\frac{Z_{\q}(x, \F_{\q + \p})}{Z_{\q}(b, \F_{\q + \p})}Z_{\q}(b, \th)\right)\p ( \q+ \p -  \k(\th))^{-1}.
\end{equation}\fn[4]{Note that in \cite{LRZ-0,LRZ,BPPR},    the function
$Z_{\q}(x,\Phi_{\q+\p})$ is denoted by
 ${H}_{\q,r}(x) = \mathrm e^{\Phi_{\q+\p}x} \left( 1 - r \int_0^{x}  \mathrm e^{-\Phi_{\q+\p}y} W_{\q}(y) \mathrm dy \right)
$}

When $\p \to \I,$ this recovers \eqref{sevruin},
and when $b \to \I$ it yields  \cite[(14)]{AIZ}
\begin{eqnarray*}&&\E_x [e^{\th X(T_0^-)} ; T_0^- < \I]=\le( Z(x,\th) - Z(x, \F_\p) \fr{ (\F_{\p+\q} -\F_\q) \k(\th)}{\p(\th -\F_\q)}\ri) \p ( \q+ \p -  \k(\th))^{-1},\end{eqnarray*}by using { \beq\label{mis} Z_\q(b, \F_{\p+\q})^{-1} {Z_\q(b,\th)} \mathop{\to}_{b \to \I} \fr{ (\F_{\p+\q} -\F_\q) \k(\th)}{\p(\th -\F_\q)}.\eeq}

\im {\bf The joint law of  dividends over an upper barrier, of the ruin time, and of the severity of ruin}. Let $E^b_x$ denote the law of $X$ reflected from above at $b$, and let $R_b(t)$ denote the regulator.  From \cite[Thm 6]{IP},
we find
\begin{eqnarray}&&\E^b_x [e^{-\vartheta R_b(\t_0^-) + \th X(\t_0^-)} ; \tau_0^- < \kil_\d]=\E _{Y(0)=b-x} [ e^{- \vt R_b(t_b) -\th (Y(t_b)-b)}] \no\\&&=Z(x,\th) - W(x) \le(W'(b_+)+ \vartheta W(b)\ri) ^{-1}  \le( Z'(x,\th)+ \vartheta Z(b,\th)  \ri) \label{e:divtill0}\end{eqnarray}

\begin{remark}  As a check, decompose in the two cases $\t_0^- < \t_b^+$, $\t_0^- \geq  \t_b^+$:
\begin{eqnarray*}&&\E^b_x [e^{-\vartheta R_b(\t_0^-) + \th X(\t_0^-)} ; \tau_0^- < \kil_\q ]=Z(x,\th) - W(x) W(b)^{-1}Z(b,\th) +\\&& W(x) W(b)^{-1} \le( Z(b,\th) - W(b) \le(W_+'(b)+ \vartheta W(b)\ri) ^{-1}  \le( Z'(x,\th)+ \vartheta Z(b,\th)  \ri)\ri)
\\&& =Z(x,\th) - W(x)  \le(W_+'(b)+ \vartheta W(b)\ri) ^{-1}  \le( Z'(x,\th)+ \vartheta Z(b,\th)  \ri)\end{eqnarray*}
\end{remark}

With $\th=0$, $x=b,$ this yields the {\bf Laplace transform  of the total dividends before ruin}
\begin{eqnarray*}&&\E^b_b [e^{ - \vt R_b(\t_0^-)} ; \tau_0^- < \kil_\q ]=Z(b) + \le(W'(b_+) W(b)^{-1}+ \vartheta I\ri) ^{-1}  \le(-Z'(b) -  \vartheta Z(b)  \ri)= \nonumber \\&&\le(W'(b_+) W(b)^{-1}+ \vartheta I\ri) ^{-1}  \le(W'(b_+) W(b)^{-1} Z(b)-Z'(b)  \ri)= \nonumber \\&&\le(W'(b_+) W(b)^{-1}+ \vartheta I\ri) ^{-1}  \le[ (W'(b_+) W(b)^{-1} - \le (W'(b_+) W(b)^{-1} \ovl W(b)- W(b)\ri )\K(0)  \ri], \nonumber \end{eqnarray*} which is the Laplace transform of a matrix exponential law. For Markov modulated L\'evy processes without discounting, $\d_i = 0, \forall i \Lra \K(0)=0$, and this reduces further to \beq\label{e:divexp} \E^b_b [e^{ - \vt R_b(\t_0^-)} ; \tau_0^- < \I ]=\le(W'(b_+) W(b)^{-1}+ \vartheta I\ri) ^{-1} W'(b_+) W(b)^{-1}.\eeq

Thus,  $R_b(\tau_{0}^-)$ has a matrix exponential distribution
generated by $\L(b)=-W_+'(b)W(b)^{-1}$.

When $\vt=0$ and with equal discounting $\d_i = \d, \forall i$, \eqref{e:divtill0} yields the {\bf severity of ruin for a regulated process}:
\begin{eqnarray}&&\E^b_x [e^{ \th X(\t_0^-)} ; \tau_0^- < \kil_\q ]=Z(x,\th) -  W(x) W'(b_+) ^{-1}  Z'(b,\th).  \label{regsevruin} \end{eqnarray}

This formula is also called the {\bf dividends-penalty} identity \cite{gerber2006note}.

If furthermore $\th=0$, this yields the {\bf Laplace transform of the ruin time  for a De Finetti regulated process},  due in  the Levy case to  \cite{AKP}:
\begin{equation}  \label{e:regruin} \P^b_x [ \tau_0^- < \kil_\q ]=Z(x) -   W(x) W'(b_+) ^{-1}   Z'(b).\end{equation}

{}

 {A further generalization in the Levy case with  Poissonian observed ruin time $T_0^-$} is  \cite[(23)]{AIZ}:
\begin{eqnarray*}&&\E^b _{x} [ e^{- \vt R(T_0^-) +\th X(T_0^-)}; T_0^- < \kil_q]=( 1- \n^{-1} \k(\th))^{-1}\\&&\le(Z(x,\th) - Z(x,\F_\n) ((\F_\n+\vt) Z(b,\F_\n) - \p W(b))^{-1}( Z'(x,\th) + \vt  Z(b,\th) )\ri).\end{eqnarray*}


\im{\bf The Gerber-Shiu function common to the  two-sided absorbed and reflected severity of ruin problems}.  Note the similarity
between the equations \eqref{sevruin}, \eqref{regsevruin} for the Laplace transforms of the severity of ruin before seeing $b$, and with reflection at $b$,  which involve the same functions $W(x), Z(x,\th)$.

\cite{APP15} show in the SNLevy case that this continues to be the case for any  pay-off $w$ which is "admissible" (satisfies certain integrability condition).
\begin{Pro}\label{prop:two}
Given  $a<b<\infty$,  $x\in(a,b)$,
and an  {\em  admissible pay-off} $w:(-\infty,a]\to\mbb R$,
there exists a unique "smooth GS function" $F_{w}$ so that the following hold:
\beq&&
V_{w}(x)=\E_x\le[e^{-q
\t_a^-}w\le(X_{|a,b|} (\t_a^-)\ri)\mbf 1_{\{\t_a^-<\t_b^+\}}\ri]  = F_{w}(x-a) - W^{(q)}(x-a)
\frac{F_{w}(b-a) }{W^{(q)}(b-a)} , \no \\&&
\T V_w(x)=\E_x\le[e^{-q
\t_a^-}w\le(X_{|a,b]} (\t_a^-))\ri)\ri]  = F_{w}(x-a)-
W^{(q)}(x-a) \frac{F_{w}'(b-a) }{W^{(q)\prime}(b-a)}.\label{e:GS}
\eeq
\end{Pro}

   Stated informally, this amounts to the fact that both these  problems admit decompositions involving an identical "non homogeneous  solution" $F_{w}$.

 The   smooth GS function
 \bea Z(x,\th) =\sum_{k=0}^\I \fr{\th^k}{k!}   Z_k(x)  \eea  corresponding to $w(x)=e^{\th x}, x \leq 0$ may be used
  as a generating function for  finding GS functions
  $ Z_k(x) $
associated to power payoffs $w(x)=x^k$.
Especially interesting are the cases $k=0,1$ which intervene in problems with linear bailout costs. Using $Z'(x,\th)=\th Z - W(x) \K(\th)$,  we find:
\beq \label{e:derlincost} && Z_0(x):= Z(x,0)= Z(x),   \;
Z_1(x):= \fr{\partial}  {\partial \th } Z(x,\th)_{\{\th=0\}}=
\ovl Z(x) -\ovl W(x) \K'(0).\eeq

The derivatives with respect to $x$ are :
\bea &&  Z'_0(x)= \q  W_\q(  x), \; Z'_1(x)= \frac{ \partial^2} {\partial x \partial \th } Z_\q(  x,\th)_{\{\th=0\}}=  Z(x) -W(x) \K'(0). \eea
\bigskip

\im {\bf Further results when   first passage is only monitored  at Poisson
  times}.
 Elegant generalizations of the classic first passage results when   observations are
 made at Poisson
  times with rate $\p$ are given in
 \cite{AIrisk,AIpower,AIZ,albrecher2015strikingly} (in the Levy case mostly):
 \BEN

\im {\bf Severity of ruin with additional killing above the barrier}. Cf. \cite[(18), Rem.2]{AIZ}, in the Levy case it holds that:
 \beq
&&\label{eid} \E_x [e^{\th X(\t_0^-)} ; \t_0^- < T_b^+ ]=\E_x [e^{\th X(\t_0^-)- \p \int_0^{\t_0^-} 1_{X(t) >b}} ; \t_0^- < \I ]\\&& Z(x,\th) - W(x) \le(\fr{\k(\th)}{\th- \F_\p}-\fr{\p Z(b, \F_\p)^{-1} Z(b,\th)}{\th- \F_\p}\ri)\no.\eeq

When $\p \to \I,$ this recovers \eqref{sevruin}, 
and when $b \to \I$ it recovers \eqref{e:old} by \eqref{mis}.

\im The solution of  the non-smooth two-sided exit problem,  when both exit times are monitored at Poisson times with equal frequency, is \cite[(19)]{AIZ}:
\begin{eqnarray*}&&\E_x [e^{\th X(T_0^-)} ; T_0^- < T_b^+]=\le(Z(x,\th) - Z(x, \F_\p) \fr{\H W_b[\F_\p, \th]} {\H W_b[\F_\p, \F_\p]}\ri)\k(\th) (\p I - \k(\th))^{-1},\end{eqnarray*}where $f[x,y]$ denotes a Newton divided difference.

\EEN

\EEN

%
{\bf Acknowledgements:}  We  thank H. Albrecher, B. Avanzi, L. Breuer, E. Frostig, J. Ivanovs, R. Loeffen,  M.Pistorius and T. Rolski  for  useful discussions.

\footnotesize

\bibliographystyle{alpha}
\bibliography{Pare31}

\end{document}